# Nonparametric Maximum Entropy Probability Density Estimation

Jenny Farmer[1,2] and Donald J. Jacobs[1,3]

[1]Department of Physics and Optical Science, [2]Department of Bioinformatics and Genomics
and [3]Center for Biomedical Engineering and Science
University of North Carolina at Charlotte, Charlotte, NC 28223, USA.

**Abstract:**

Given a sample of independent and identically distributed random variables $\{V_k\}$, a novel nonparametric maximum entropy method is presented to estimate the underlying continuous univariate probability density function (pdf). Estimates are found by maximizing a log-likelihood function based on single order statistics applied to $\{U_k\}$ after transforming $V_k \to U_k$ through a sequence of trial cumulative distribution functions that iteratively improve using a Monte Carlo random search method. Improvement is quantified by assessing the $\{U_k\}$ random variables against the statistical properties of sampled uniform random data. Quality is determined using an empirically derived scoring function that is scaled to be sample size invariant. The scoring function identifies atypical fluctuations, for which threshold values are set to define objective criteria that prevent under-fitting as trial iterations continue to improve the model pdf, and, stopping the iteration cycle before over-fitting occurs. No prior knowledge about the data is required, although boundary conditions, symmetry constraints and censor windowing can be optionally specified. An ensemble of pdf models is used to reflect uncertainties due to statistical fluctuations in random samples, and the quality of the estimates is visualized using scaled residual quantile plots that show deviations from size-invariant statistics. These considerations result in a tractable method that holistically employs key principles of random variables and their statistical properties (e.g. maximum entropy, order statistics, log-likelihood estimate, sample size fluctuations) combined with employing orthogonal basis functions and data-driven adaptive algorithms. Benchmark tests show that the pdf estimates readily converge to the true pdf as sample size increases. Robust results are demonstrated on several test probability densities that include cases with discontinuities, multi-resolution scales, heavy tails and singularities in the pdf, suggesting a generally applicable approach for statistical inference.

# 1. INTRODUCTION

A central tenant of statistical modeling is that different samples of a random process yield fluctuating outcomes from a population that is characterized mathematically by a probability density function (pdf). Conversely, finding a pdf that is consistent with sampled univariate data is a common task in many applications [1-8]. In this paper, we present a nonparametric method to solve the inverse problem of estimating a pdf that appropriately describes a dataset of $N$ random variables, $\{V_k\}$, which are assumed to be independent and identically distributed (i.i.d.). No other information about the input data is needed, although symmetry constraints and boundary conditions can be specified. As an inverse problem, the estimated pdf from a finite dataset does not have a unique solution, although variations among solutions decrease as sample size increases. Therefore, the pdf that we give as an estimate is selected from an ensemble of possibilities, where uncertainty in the estimate is quantified by a scoring function that is sample size invariant, and, to our knowledge, is introduced for the first time. The method resists over fitting to sample fluctuations and it resists under fitting the statistically significant features in the pdf caused by over smoothing. A simple analytical expression for the estimated pdf allows for subsequent calculations of statistical properties. Without any human subjective judgments, the method performs markedly well on several probability distributions with disparate difficult characteristics, where the quality of the estimates is benchmarked using four distinct metrics. In all cases, the pdf estimates improve as sample size increases. This work establishes a new route for probability density estimation that is fully automated and data driven, both being critical features sought in high throuhput applications [9-12].

## 1.1 Methods for Density Estimation

Estimating a probability density function (pdf) from random sampled data is a well-studied problem in statistics [1, 3, 13, 14]. The simple method of binning the data to obtain a histogram offers insight into how the data is distributed, but this approach falls short of producing an analytical function to represent the probability density. A problem encountered with binning data is in the choice of the bin width, which affects the appearance of the histogram [15-20]. Too small of a bin size makes the histogram noisy because the small number of samples that fall within a bin creates relatively large fluctuations. On the other hand, too large of a bin size will result in an overly coarse-description, losing characteristic features about the data. In applications where the data derives from a single resolution scale, a single bin width is sufficient to obtain insightful results. However, the selection of bin width becomes challenging when the data exhibits multiple characteristic resolution scales (i.e. such as a mixture model of Gaussians with widely varying standard deviations) where the advantages and disadvantages of choosing small and large bin widths conflict. Hence a data-driven adaptive bin width strategy is required to balance these objectives.

It is useful to replace a histogram with a smooth function representing a pdf for the data to reflect a consensus of the data points that fall within local regions. A common approach to this is kernel density estimation (KDE), which constructs a smooth curve that tracks the shape of a histogram [21] using a linear combination of kernel functions. The analytical function from KDE serves as a *possible* probability density function of a population from which the sampled data was drawn. A bandwidth associated with the kernel function controls how much data is clumped together locally. Among a vast number of potential kernel functions, the Gaussian is popular, and often a single bandwidth is used. Distinguishing between fluctuations in a random sample versus the result of features in the true pdf is challenging. If multiple resolution scales are needed, it is difficult to smooth out unwanted noise while retaining sharp features not caused by fluctuations. Although good smoothing methods are available to account for multiple resolution scales [22, 23], subtle problems inherent to KDE include choosing the kernel [24], dealing with boundary conditions [25], and objectively determining multiple resolution scales [26, 27].

The maximum entropy method (MEM) [28-30] is another approach to estimate the analytical form of a pdf described by $p(v)$ on support $[a, b]$. The unknown function $p(v)$ is determined using calculus of variations [31] by maximizing the entropy, $S = -\int_a^b \ln(p(v)) p(v) \, dv$, under the normalization constraint given by $\int_a^b p(v) \, dv = 1$. This normalization constraint restricts the allowed solution space to integrable functions. For this minimally constrained problem, the solution is $p(v) = \frac{1}{b-a}$, meaning the pdf is uniform on $[a, b]$. If the only information known about the random sampled data is that it falls within $[a, b]$, then the pdf that maximizes entropy assigns equal probability to every equally spaced differential element within $[a, b]$, and hence $p(v)$ is a constant. If certain moments are known, this information is captured as additional constraints, which causes $p(v)$ to change shape as appropriate [32]. Many different types of constraints can be employed, such as a set of power moments. For example, if only the first and second moments of the random variable, $V$, are constrained, it is well known that the shape of the pdf becomes that of a Gaussian distribution. More generally, the shape of the pdf can be systematically modified by including higher order constraints through a complete set of orthogonal basis functions [33].

**1.2 Motivations and Challenges**

Regardless of the approach employed, any pdf that is inferred from a random sample of data will have imperfections due to incomplete knowledge extracted from a finite sample, and inexactness from the numerical method employed [34]. Our original motivation to develop a nonparametric pdf based on maximum entropy was spawned from the desire to accurately calculate the density of states as a function of energy in large molecular systems [35-37]. The problem of calculating the density of states was transformed to finding an accurate probability density. To that end, a robust MEM [33] was developed to model this probability density, especially within low energy tails that are difficult to determine. While successful in providing high accuracy for sample sizes with 100,000 or more observations in the original application, the method was not resistant to the general problem that MEM has in over fitting data [38, 39], especially for small sample sizes. Although heuristic criteria were developed to avoid over fitting, no objective metric was established to lift subjectivity completely. Significant improvement for establishing an objective approach that eliminates the need for subjective human intervention is reported here. However, before discussing these new developments, it is instructive to first describe the origin of the problems that challenge the framework of any MEM.

The MEM potentially requires knowing high order moments that may well exceed $\langle v^{100} \rangle$. High order moments often do not exist, and when they do exist, they tend to cause calculations to become unstable due to amplification of statistical errors when moments are empirically estimated [40]. The source of this error derives from sample fluctuations where an empirical average differs from the population average. This is troubling for the inverse problem because moments for the true pdf are not known. Rather, moments are estimated based on sample averages containing statistical errors. Therefore, when large power-moments are needed, the problem often becomes ill-conditioned and/or numerically intractable. The immediate convergence problems caused by power-moments are avoided by considering moments of level-functions [33, 36]. Level functions have lower and upper bounds. As an example, $\sin(av)$ and $\cos(bv)$ are level functions because they are bound between -1 and 1. The moments of level-functions always exist, which reduces numerical instability.

Another problem with a traditional MEM is that it is formulated as a parametric model [41], where the number of moment constraints is known in advance. Usually constraints are defined in terms of a set of functions, $\{g_j(v)\}$, which are meaningful for the specific problem being solved. As briefly reviewed in Section 2.1, MEM leads to $exp[\sum \lambda_j g_j(v)]$ as the canonical form for a pdf, where $\{\lambda_j\}$ is a set of series coefficients that are determined in a data driven way. Of course selecting constraint functions based on knowledge of the underlying process is not possible if nothing is known about the data. Instead, a set of level-functions, $\{g_j(v)\}$, will be employed that also define a complete basis set of orthogonal functions on the domain of interest. This ensures that the series expansion $\sum \lambda_j g_j(v)$ can model any continuous integrable function. It is worth mentioning that although orthogonal functions have been used in KDE [42], the expansion within MEM appears in the argument of an exponential function as a consequence of satisfying the maximum entropy condition. As higher order level-functions are included in the series expansion, sharper features in the pdf can be captured because higher order level-functions oscillate more rapidly. Conceptually, it is easy to see that sample size places a maximum limit on the number of terms that can be justified in the expansion. Specifically, rapidly oscillating functions will at some point characterize features on scales that are much finer than the data itself can support. Although this limit exists, determining how many level-functions to include in the series expansion is a critical problem that is solved in this paper.

In practice, fluctuations in a sample will influence the number of terms needed to accurately model the data. However, attempting to model a particular sample well is likely to lead to over-fitting because statistical fluctuations will be captured that are not representative of the population [5, 43, 44]. As such, higher orders of level-functions are justified commensurately as more data is collected, and this provides the statistical resolution necessary to estimate sharper features in a pdf [45]. Unfortunately, working with a finite series expansion for $\sum \lambda_j g_j(v)$ also generally causes wiggles to appear in the model pdf that are not present in the true pdf. Although the condition to maximize entropy suppresses wiggles that appear in model solutions, it is possible to include smoothing constraints within the framework of a MEM [33, 46]. However, imposing smoothing criteria is a subjective constraint that is not data driven.

The nonparametric MEM reported here differs from an earlier working version [33] in substantial ways involving fundamental principles of statistical science. However, there are two aspects that remain similar in both methods. First, Chebyshev polynomials are selected for the level-functions employed here, and constitute the orthogonal basis set. Second, the expansion coefficients are determined by employing a random search method [47-49] to optimize a multivalued non-linear scoring function. The random search method, called funnel diffusion, has been employed and explained previously [33]. All other details about the non-parametric MEM presented here are quite different from our earlier work.

**1.3 New Contributions**

Given a set of i.i.d. random variables $\{V_k\}$, their pdf is estimated by iteratively adjusting expansion coefficients to increase a log-likelihood function [50, 51] for single order statistics of the transformed random variables $\{U_k\}$. Here, $\{U_k\}$ is isomorphic to $\{V_k\}$ by transforming $V_k \rightarrow U_k$ using a trial cumulative distribution function [52] associated with a trial $p(v)$. This mapping always exists, and it is found to be numerically stable. Integrating $p(v)$ is done efficiently through a data-driven adaptive algorithm. An important advance is in constructing a sample size invariant scoring function to quantify how good a trial pdf is. This scoring function is used to monitor how well $\{U_k\}$ represents *sampled uniform random data* (SURD). Furthermore, an efficient algorithm results by taking advantage of single order statistics across uniformly spaced subsets of sampled data [53]. Consequently, hierarchical partitions of sort ordered datasets (typically divided by factors of 2) is considered, starting with the smallest partition (i.e. sparse data) and ending with all data in the last partition analyzed.

From this data-driven iterative Monte Carlo procedure, an ensemble of possible model solutions is generated, from which an estimate is made, and uncertainties can be readily quantified. Although uncertainties in KDE estimates can be quantified [54], in practice it is usually overlooked and under-appreciated. In this work, our focus is to test the quality of all members in the ensemble. Therefore, scaled quantile residual plots (SQR-plots) are introduced to visualize the quality of a model pdf. Unlike a standard QQ-plot that approaches a straight line with slope of 1 going from 0 to 1 on both the abscissa and ordinate as sample size increases [55], a SQR-plot is sample size invariant. Finally, a figure of merit is introduced to provide another metric to quantify the goodness of a pdf estimate.

**1.4 Remaining Layout of Paper**

In Section 2, preliminary steps for developing a pdf estimate are described that include a description of MEM, how the domain for the pdf is determined, and how functional space of continuous integrable functions is explored. Interestingly, moments of level-functions are not calculated in favor of iteratively calculating a trial cdf. In Section 3, data-driven adaptive methods to calculate the cdf from the pdf are discussed, as well as single order statistics and the process employed to arrive at a novel universal scoring function based on a log-likelihood function for SURD. In Section 4, the numerical methods employed to facilitate the calculations are briefly discussed. In Section 5, results for eight example distributions that embody difficult facets are discussed. In Section 6, conclusions are made based on the numerical results that demonstrate the nonparametric MEM is robust, versatile, and consistent with statistical resolution and fluctuations expected in finite random samples of data.

## 2. NONPARAMETRIC STATISTICAL MODEL

**2.1 Maximum Entropy Probability Density**

Calculus of variation is employed to determine an appropriate form for a pdf that maximizes entropy subject to a variety of constraints using the method of Lagrange multipliers. For a general approach, the constraints are expressed in terms of moments of a set of functions $\{g_j(v)\}$ comprising a complete set of orthogonal level-functions. In principle, moments of these level-functions can be calculated for any j-index, but at this stage it is not a concern how $\langle g_j(v) \rangle$ is obtained. Rather, only one question is being asked: What is the general form of $p(v)$ that is consistent with all constraints while maximizing entropy? Viewed as a function of $p(v)$, the functional $\tilde{S}$ comprised of the entropy, $S$, with constraint conditions added, is explicitly written as:

$$\tilde{S}[p(v)] = \lambda_0 \left( \int p(v)\, dv - 1 \right) + \sum_{j=1}^{D} \lambda_j \left( \int g_j(v)\, p(v)\, dv - \langle g_j \rangle \right) - \int \ln[p(v)]\, p(v)\, dv \tag{1}$$

The variables, $\lambda_j$, are Lagrange multipliers that play an important role in determining the shape of $p(v)$ based on knowing the exact population averages $\{\langle g_j \rangle\}$ for which the $\lambda_j$ are conjugate to. In the process of maximizing entropy, the $-\ln[p(v)]\, p(v)$ term in the integrand constrains all viable $p(v)$ functions to be non-negative. The solution to Eq. (1) for $p(v)$ is straightforward to obtain [33], which works out to be

$$p(v) = \exp\left[ (\lambda_0 - 1) + \sum_{j=1}^{D} \lambda_j g_j(v) \right]. \tag{2}$$

In Eq. (2) the Lagrange multipliers are unknown variables that must be determined (somehow). It is worth pointing out that solving for the Lagrange multipliers is the difficult part of the process, but once solved, the values of the Lagrange multipliers provide insight into the importance of a constraint. For instance, $\lambda_j = 0$ indicates $\langle g_j \rangle$ is an unnecessary constraint. Larger $|\lambda_j|$ indicates a greater affect on how the shape of the pdf deviates from uniform. Our intuition views these Lagrange multipliers as constraint forces in analogy with classical mechanics. Regardless of the interpretation, the problem at hand is to solve for them. In addition, the number of Lagrange multipliers is also unknown.

The general form of the pdf given in Eq. (2) shows that in all cases, $\lambda_0$ ensures the total probability will be 1. Consequently, $\lambda_0$ is dependent on $\lambda_1$ to $\lambda_D$ where $D$ defines the dimension of a parameter space. Interestingly, $\lambda_0$ is the only Lagrange multiplier that scales the pdf while preserving its shape. It is convenient to define $\Lambda = \lambda_0 - 1$, such that $e^\Lambda$ is a normalization constant, and to let $\boldsymbol{\lambda} = (\lambda_1, \cdots \lambda_D)$ as a simpler notation. Here, $\boldsymbol{\lambda}$ is a point in parameter space of $D$ dimensions, having $\Lambda(\boldsymbol{\lambda})$. Furthermore, the solution from MEM is now written as $p(v|\boldsymbol{\lambda})$ to emphasize the model pdf is a function of $D$ independent parameters. An important observation is that a finite $D$ implies details on small scales are unconstrained by the model pdf since $\lambda_j = 0$ for $j > D$. To justify working with a larger $D$ (i.e. more parameters need to be determined) requires more data in order to have the necessary statistical resolution to impose a constraint. The connections between statistical resolution, capturing the fine structure of a pdf, and a data-driven process that automatically justifies increasing the dimension of parameter space, derives in large part from the mathematical properties of expanding in a complete set of orthogonal functions.

**2.2 Bound Intervals and Censor Windowing**

When working with numerical estimates, describing a pdf on a bound interval $[a, b]$ does not pose much limitation in practice. First, if bounds on the random variables are known, they should be used to obtain maximum accuracy for the model pdf. Second, selecting a finite domain is necessary in applications where samples are censored by windowing. For example, a selection rule may measure the value of a random variable $V_k$ within the limits $[a, b]$, while counting the number of events when $V_k$ falls outside this range. The ratio, $R_{ab}$, of the number of events measured within the range $[a, b]$ to the total number of events gives an estimate for the probability, $P_{ab} = \int_a^b p(v)\, dv$. The normalization constraint required in the MEM generalizes to $P_{ab} = R_{ab}$, where uncertainty is introduced from the empirical estimate for $P_{ab}$. Fortunately, the precise location of the interval $[a, b]$ is not critical because a model pdf can be accurately calculated within the censoring window regardless of how much data is discarded, provided the amount of discarded data is counted. Although not explored in this work, it is worth noting that the concept of censor windowing could, in principle, be used to develop a parallelized version of estimating a pdf by applying a series of adjacent censoring windows in sequence, where each window is processed independently, and at the end the pieces are stitched back together.

In applications without a priori knowledge of the limits on $V_k$, it generally suffices to extend the range beyond the observations such that $a < \min\{V_k\}$ and $b > \max\{V_k\}$ because this extended range helps justify that the probability outside the range is negligible. We use a data driven criteria for this. For a sample of N data items, sorted from lowest to highest values, we set $a = V^{(1)} - [V^{(5)} - V^{(1)}]$ and $b = V^{(N)} + [V^{(N)} - V^{(N-5)}]$. Here, the sort index is used as a superscript. In words, the left lower limit is shifted further to the left of the lowest observed data item by the difference between the 5[th] lowest and lowest observed values in the dataset. Similarly, the right upper limit is shifted further to the right of the highest observed data item by the difference between the highest and 5[th] highest observed values in the dataset. However, when the true pdf has a heavy tail describing extreme statistics, the range must be reduced to obtain stable numerical integrations. This means, $\min\{V_k\} < a$ and/or $\max\{V_k\} > b$ so that $[a, b]$ defines a censoring window. A balance must be reached to cut off a low probability tail to ensure enough data samples are present to make an accurate estimate, while not throwing out too much data in the tail, which would otherwise lose information about the tail. To balance these disparate objectives, the MEM presented here considers a closed interval $[a, b]$ that keeps nearly all data, except for outliers according to the formulas: $a = Q_{25} - c(Q_{75} - Q_{25})$ and $b = Q_{75} + c(Q_{75} - Q_{25})$ and $Q_{25}$ and $Q_{75}$ are the 25% and 75% quantiles respectively. When $c = 3$, the outliers identified have been called *extreme* [56, 57]. For this work, *very extreme* outliers are identified using the default setting of $c = 7$. If no extreme outliers exist (i.e. no heavy tails) all data is considered and the limits on the interval are larger than the random data points span. For a heavy tail, extreme outliers must be considered to capture the important feature of the true pdf. This default cutoff works well because extreme statistics is captured with enough sample statistics to infer a good estimate for the probability in the tails, and, as sample size increases, the interval $[a,b]$ being considered is likely to stretch outward using this data driven quantile approach.

Without loss of generality, a linear transformation is employed to map the random variable $V_k$ to $X_k$ according to $X_k = \frac{2V_k - b - a}{b-a}$. This transformation conveniently maps the domain onto the interval $[-1,1]$, which is natural for Chebyshev polynomials [58]. In particular, Chebyshev polynomials of the first kind are level-functions, because $|T_j(x)| \leq 1$ for $-1 \leq x \leq 1$ and for all $j$. Moreover, the set of functions, $\{T_j(x)\}$, form a complete orthogonal basis on $[-1,1]$. Although Chebyshev polynomials have been selected because of their convenient properties, other sets of level-functions that form a complete set of orthogonal functions on a bounded interval could be used. As the last step in the process, the solution for $p(x)$ is transformed back to the original variable. Unless stated otherwise, henceforth the discussions of the method will be with respect to the transformed variables $X_k$ and their pdf, denoted simply by $p(x)$.

**2.3 Infinite Dimensional Function Space**

A model pdf will be denoted as $p_e(x)$ because it is an estimate for $p(x)$. In more explicit detail the dependence on parameter space is denoted as $p_e(x|\lambda)$ to highlight the idea that the shape of a model pdf changes as the point $\lambda$ moves through parameter space. In principle, the dimension of parameter space may have to be infinite to represent the true population pdf. To see why, define $h(x) = \ln[p(x)]$ on the interval $[-1,1]$. It follows that $h(x) = \sum \lambda_j T_j(x)$. The exact values of $\{\lambda_j\}$ can readily be determined uniquely by spectral decomposition. Since $\{T_j(x)\}$ is a complete set, the form of Eq. (2) should always provide an adequate representation of the true pdf. However, a problem occurs for values of $x$ where $p(x) = 0$. In this case, $|h(x)| \to \infty$, and if $p(x) = 0$ for any finite stretch of $x$ within $[-1,1]$, then $h(x)$ is not a square-integrable function. However, in practice $p_e(x)$ will not equal zero when $p(x) = 0$. Instead, if an acceptable tolerance is $\epsilon$ (say $\epsilon = 10^{-12}$), then a pragmatic answer will be achieved by a series expansion when it is equated to $h_e(x) = \ln[p(x) + \epsilon]$. This and other numerically bad cases, such as when $p(x)$ has a singularity, are handled by setting tolerances for the desired numerical accuracy.

As a consequence of Eq. (2), $p_e(x)$ is a continuous function that is infinitely differentiable, integrable, and importantly, the form of Eq. (2) represents a general approach. Specifically, the series expansion in Chebyshev polynomials for $h_e(x)$ will yield an estimate of $p(x)$ through $\exp[h_e(x)]$ that will be accurate to within some numerical precision. This precision is limited by statistical resolution set by sample size. For a finite sample, a finite number of Lagrange multipliers will be needed. A data-driven procedure must be developed to objectively determine the dimension of parameter space, as well as the Lagrange multipliers corresponding to a model pdf that would be likely to have generated the sample.

**2.4 Determination of Lagrange Multipliers**

In previous work [33], Lagrange multipliers defining a point in parameter space were determined by funnel diffusion, which is a simple adaptive random search method (see Appendix A.1). Starting at $\lambda^{(0)}$ the algorithm makes trial random steps in the parameter space that land at a new $\lambda$. A trial step produces $p_e(x|\lambda)$ that may have a functional form that is farther or closer to $p(x)$ as quantitatively measured by a scoring function. The trial step is taken when the score improves; otherwise the trial step is not taken. Labeling only the successful steps, $\lambda^{(i)}$ indicates that $p_e(x|\lambda^{(i)})$ is closer to $p(x)$ than $p_e(x|\lambda^{(i-1)})$ based on a scoring function. Previously [33], a least squares error (LSE) was calculated between moments based on the current pdf estimate and the empirical moments using a plugin method. Although the Monte Carlo approach using funnel diffusion minimizes the LSE efficiently, as the sample size decreases, uncertainties in the target moments increase, creating deviation from $p(x)$ due to statistical fluctuations inherent in a small dataset. Consequently, the previously employed MEM tends to over fit the data because the target moments reflect fluctuations in the dataset itself.

Within a maximum entropy framework, the problem of over fitting would not be present if the empirically estimated target moments could be replaced with actual population moments, or perhaps some other statistical metrics determined by the population. The fundamental improvement described in section 3 for the MEM is that the error between a trial pdf and the true pdf can be calculated using a universal statistical metric without requiring knowing $p(x)$. The word "universal" is used to emphasize that this scoring function is problem independent. It is also approximately independent of sample size. The method of funnel diffusion yields an ensemble of possible solutions that are close to the exact solution while accounting for expected levels of fluctuations in finite samples. As the size of the sample increases, the possible solutions in the ensemble are found to converge to the true pdf. Another improvement employed is that the number of Lagrange multipliers (parameters) is not a priori specified.

**3. UNIVERSAL STATISTICS BENCHMARK**

Statistics on *sampled uniform random data* (SURD) on the interval $[0,1]$ will be used as a universal benchmark to characterize how good a trial pdf is at describing the random sample of data. The strategy employed is to transform the i.i.d. random variables $\{X_k\}$ that span the bounded interval $[-1,1]$ onto new random variables $\{U_k\}$ that span the interval $[0,1]$ using the cdf that is associated with a trial pdf. A trial pdf will be considered an appropriate approximation to the true pdf when $\{U_k\}$ exhibits the same statistical properties as that of SURD. Therefore, quantitative metrics are constructed that compare certain statistical properties of $\{U_k\}$ to SURD. Several types of statistical properties were tested that included various moments, binning statistics, and order statistics.

Binning statistics is based on uniformly dividing the interval from 0 to 1, into a fixed number of bins, where each bin on average will hold approximately 10 sample points. A binomial distribution gives the probability for a specific number of data points to fall within a bin. Details of how single order statistics is employed are explained in the next section. Interestingly, single order statistics was found to be most effective in distinguishing between typical and atypical fluctuations in sample datasets. It is worth noting that tracking multiple statistical properties simultaneously, such as combining binning statistics with single and/or double order statistics, does not offer a clear advantage across several diverse test-distributions. Therefore, the simplest approach of using single order statistics is presented in this report.

As explained in section 2.4, funnel diffusion anneals the estimate $p_e(x|\lambda)$ toward the true pdf as the parameter space is explored, guided by the scoring function. A funnel shape landscape characterizes the scoring function as the $\lambda$ parameters deviates from the exact values (but unknown). Conceptually, the bottom of a funnel represents potential solutions consistent with fluctuations from a finite sample. This funnel sharpens as sample size increases, and for small sample size, a more blunt or rounded funnel occurs because it reflects more typical fluctuations. Operationally, the cdf is determined from the pdf, as $F_e(x|\lambda) = \int_{-1}^{x} p_e(t|\lambda)\, dt$. The cdf is used to transform the random variables $\{X_k\}$ to a new set of isomorphic random variables $\{U_k\}$, where $U_k = F_e(X_k|\lambda)$. Regardless of the Lagrange multipliers at the current step in the annealing process, it is guaranteed that all $\{U_k\}$ will be in the interval [0,1]. If $p_e(x|\lambda)$ is the true pdf, then $\{U_k\}$ will be SURD on [0,1]. Conversely, if the statistical properties of $\{U_k\}$ are very different than what is expected, continuation of generating and checking new trials is necessary to improve the trial pdf. The goal is to iteratively adjust the Lagrange multipliers to improve the match between order statistics of $\{U_k\}$ and SURD. The process is terminated when the match is within expected levels of fluctuations, which has a simple dependence on sample size.

**3.1 Order Statistics and Scaled Residual Quantile Plots**

Order statistics provide a sensitive metric for characterizing SURD on [0,1]. The random data $\{U_k\}$ assumed to represent SURD on [0,1] is sorted from smallest to largest values. The sorted data is labeled as $\{U^{(s)}\}$ where $s$ is a sorting index from 1 to $N$, such that $U^{(s)} < U^{(s+1)}$ for all $s$. The properties of order statistics for a uniform probability distribution are well known in the literature [59-61]. Employing single order statistics on $N$ observations, the probability density for finding $U^{(s)}$ at $u$, is given by

$$p_s(u|N) = \frac{N!(1-u)^{N-s}u^{s-1}}{(N-s)!(s-1)!} \tag{3}$$

From Eq. (3) it follows that the mean position for $U^{(s)}$ and its variance are respectively given by

$$\mu(s|N) = \frac{s}{N+1} \qquad \sigma^2(s|N) = \frac{s(N+1-s)}{(N+2)(N+1)^2} \,. \tag{4}$$

After simplifying the variance, the standard deviation characterizing typical fluctuations in the position of $U^{(s)}$ about its mean is given by

$$\sigma = \frac{\sqrt{\mu(1-\mu)}}{\sqrt{N+2}} \,. \tag{5}$$

Note that the mean positions of single order statistics for SURD span the [0,1] interval with uniform spacing, which is an intuitive result. The fluctuation about the mean given in Eq. (5) is also interesting in that the plot of $\sigma^2$ versus $\mu$ is a concave parabola, indicating that a random variable $U^{(s)}$ near either end has much less variation due to the fixed boundaries. Importantly, the positions of the sorted random variables tend closer to predicted average locations as sample size increases, while deviations from averages decrease. In particular, typical sample size fluctuations within the set of positions $\{U^{(s)}\}$ about their respective average positions drop as $\frac{1}{\sqrt{N+2}}$ according to Eq. (5). Therefore, larger sample sizes will resolve single order statistics better, and the size of typical fluctuations provides a good way to quantify and visualize statistical resolution.

When SURD is generated in numerical experiments, deviations of the generated random variables from their theoretical population average define a residual quantile. The residual for the $s$-th sorted data point is defined by $\delta_s = U^{(s)} - \mu(s|N)$, and when plotted against the expected mean, $\mu(s|N)$, defines a residual quantile plot. Unfortunately, it quickly becomes clear that residual plots are not very useful when used to track how a model pdf converges to the true pdf as sample size increases, or to visualize the accuracy of $p_e(x|\lambda)$, until after the scale of the plot is magnified to view the residuals as they approach zero as the number of samples increase. Therefore, it is convenient and advantageous to scale the deviations by a sample size factor to arrive at a sample size invariant measure. Specifically, let

$$\Delta_s = \sqrt{N+2}\,\delta_s = \sqrt{N+2}\left(U^{(s)} - \frac{s}{N+1}\right). \tag{6}$$

The variable $\Delta_s$, when plotted against $\mu(s|N)$, defines a scaled residual quantile plot (SRQ-plot). The SRQ-plot makes it easy to visually assess the results for different sample sizes on equal footing, and systematic errors in pdf estimates can be visually identified without knowing the true pdf. Importantly, the SRQ-plot provides a model independent metric to visually characterize the quality of a pdf estimate within the expected level of statistical resolution set by sample size. As an illustrative example, SURD is generated using a uniform random number generator to create four samples of various sizes. In **Figure 1A** the fluctuations in these samples are visualized in terms of the commonly used quantile-quantile plot (QQ-plot). In **Figure 1B** the same fluctuations are visualized using SRQ-plots, which show the amplitude of the scaled fluctuations are typically of order 1 across all sample sizes.

## 3.2 Sample Size Invariant Scoring Function

To quantify the typical fluctuations that are expected in SURD, extensive numerical experiments were conducted. Using a robust uniform random number generator [62], a sample of N random numbers is generated on [0,1]. For example, the QQ-plots and SRQ-plots in **Figure 1** show particular outcomes of these numerical experiments. Each set of $N$ random numbers defines a sample, and at least 100,000 samples are made per $N$, with millions of samples generated for small $N$. The probability given by Eq. (3) for locating $U^{(s)}$ in a sample is multiplied together to define a product measure given as

$$P_L\left(\{U^{(s)}\}\right) = \prod_{s=1}^{N} p_s\left(U^{(s)}|N\right) \quad . \tag{7}$$

To simplify the calculations, the natural logarithm of Eq. (7) is taken to transform a product of terms into a sum of terms. Dividing this sum by $N$ defines an average for the entire sample. Each sample is therefore assigned a mean value that quantifies the single order statistics for the $\{U^{(s)}\}$ configuration. Numerical experiments reveal that the expectation of these mean values over an ensemble of samples is accurately described as a linear function of $\ln(N)$ such that

$$E\left[\tfrac{1}{N}\ln(P_L)\right] = \tfrac{1}{2}\left(\ln(N) - \tfrac{4}{5}\right) \quad . \tag{8}$$

By subtracting the systematic $\ln(N)$ dependence from each sample, a modified log-likelihood metric is defined as a function of a specific single order statistic configuration (i.e. the set $\{U^{(s)}\}$) as follows

$$L\left(\{U^{(s)}\}\right) = \tfrac{1}{N}\ln\left(P_L\left(\{U^{(s)}\}\right)\right) - \tfrac{1}{2}\ln(N) \quad . \tag{9}$$

The modified log-likelihood metric, $L$, turns out to be sensitive to small variations in $\lambda$. Recall that when $\lambda$ is changed the shape of the trial pdf and its associated cdf changes, which results in a change in the way $\{U^{(s)}\}$ is distributed on the interval [0,1]. The high sensitivity in $L(\{U^{(s)}\})$ on the way $\{U^{(s)}\}$ is distributed on [0,1] is a desirable property to have because it facilitates a robust assessment of a good or bad trial $p_e(x|\lambda)$. Henceforth, the word *modified* will be dropped when referring to this $L$ metric.

From Eq. (9) the log-likelihood values that come from SURD have a mean of $-0.40$ for any sample size (after subtracting out the systematic $\ln(N)$ dependence). Interestingly, the probability distribution for log-likelihood metrics across different sample sizes exhibit data collapse such that the pdf for log-likelihood is approximately sample size invariant as shown in **Figure 2**. As explained below, small variations between scaled functions do not lead to any practical concern. Therefore, a single universal scoring function is constructed from a compilation of log-likelihood metrics from different sample sizes ranging from $2^8$ to $2^{20}$ by factors of 2. In this report, the smallest sample size considered is $N = 256$ where focus is on universality of the method. While the universal scoring function is found by us to work well in practice for all $N$, application of the MEM for small sample sizes will be published elsewhere.

The maximum of the log-likelihood metric occurs when $U^{(s)}$ lands at the mode of the s-th sorted value (i.e. the $u$ where $p_s(u|N)$ is a maximum defines its mode) for all $s$. However, there is zero probability to observe $N$ random numbers on [0,1] that happen to be spread in exact accordance with all the individual mode locations. Indeed, the probability of observing this or any particular configuration of points is zero. Importantly, Eq. (9) provides a sensitive measure for characterizing how the set of points as a whole deviate from the population average of $N$ sampled uniform random numbers on [0,1]. However, with respect to this log-likelihood metric, some deviations from the population average will be more probable than others. We note that an alternative least squares error metric for the deviations in the positions of $\{U^{(s)}\}$ from the corresponding population average positions (normalized by $\sigma(s|N)$ as Z-scores) was implemented and tested. From numerical comparisons, the log-likelihood metric of Eq. (9) provides superior sensitivity for discriminating between typical and atypical $\{U^{(s)}\}$ configurations, and it affords rapid evaluation.

A particularly important issue is to model the pdf at the boundaries of the sampled data correctly. To ensure proper limits are satisfied requires augmenting a penalty term to the scoring function discussed above. Incorporating a penalty term was found to be essential for distributions with heavy tails. Although the functional form for the penalty does not appear to be critically important, all the results reported here have been obtained using a penalty that is subtracted from $L$ based on deviations from uniform at the boundaries as follows.

$$penalty = \ln\left[1 + \frac{0.1}{p}\sum_{i=1}^{p}\left|U^{(i)} - \frac{i}{N+1}\right| + \frac{0.1}{p}\sum_{j=N-p}^{N}\left|U^{(j)} - \frac{j}{N+1}\right|\right] \qquad \text{Eq. (10)}$$

Here $p = \lfloor N * 0.005 \rfloor$. This choice of $p$ means that 1% of the data (0.5% at both the left and right boundaries) is used in the penalty term. Overall, this penalty has virtually no effect on good solutions, and for these good solutions removing the penalty does not yield negative consequences. Conversely, good solutions are difficult to find if the penalty term is not augmented to the universal scoring function at the start. That is, the penalty term guides funnel diffusion very strongly when the $\lambda$ parameters are far from converged values given by the appropriate (i.e. using the selected basis functions) spectral decomposition of the natural logarithm of the true pdf. In summary, augmenting this penalty to the scoring function seeds good solutions that produce mappings consistent with SURD at the boundaries.

As shown in **Figure 2**, simulated outcomes of $L$ for SURD are accumulated to construct the pdf for the universal scoring function. Its single peak corresponds to the most probable log-likelihood that will result from SURD. Through trial and error adjustments of the Lagrange multipliers that produce different sets of $\{U^{(s)}\}$, the initial objective is to maximize the log-likelihood $L(\{U^{(s)}\})$. In general, the sensitivity of a scoring function plays a critical role in how an iterative method of any type (random search method, simulated annealing, etc.) converges toward a solution, affecting speed, accuracy and selectivity. Although it remains prudent to further explore alternatives, the log-likelihood of single order statistics is found to be markedly stable and it exhibits the best performance characteristics among several alternatives that were already tested, including least squares error, moments, and binning statistics.

### 3.3 Solutions for the Inverse Problem

As mentioned above, the main task initially is to maximize the log-likelihood function. But at some point, maximizing log-likelihood makes the $\{U^{(s)}\}$ obviously less random because it becomes atypically *perfect*. For example, suppose funnel diffusion can always find the values of the Lagrange multipliers, no matter how many are needed, such that the final values $\{U^{(s)}\}$ land at the mode positions of $p_s(u|N)$. Surely, something is wrong! Drawing a sample of $N$ random values from a known pdf will result in a set of $\{U^{(s)}\}$ that will look like SURD, meaning there will be fluctuations around the expected mean positions. Thus, *maximizing the log-likelihood is not an appropriate goal*, but rather, the objective is to set an acceptable range on the scoring function that represents typical outcomes. Indeed, because the scoring function is the probability density of log-likelihoods for a sample of $N$ data points, a confidence window can be defined. For example, to arrive at a confidence of 90%, the area under the curve for both the lower and upper tails of the scoring function can be cutoff at the 5% level. Hence, the cumulative distribution of log-likelihoods for SURD would be marked between 5% and 95%.

In practice, the process of optimizing the Lagrange multipliers has to be terminated using some criteria. If a low value of the log-likelihood is selected to stop funnel diffusion with 5% *coverage*, this means that only 5% of all log-likelihoods from SURD are less than this value. Thus, 5% coverage will produce atypical results. Similarly, atypical results are produced if the target log-likelihood is selected at 95% coverage, where 5% of all possible log-likelihoods are greater than this value. This is not to say that the pdf estimates produced at low and high target values are not possible. The problem is, why settle for only the bottom 5%, or demand only the top 5%, compared to all accessible pdf estimates that can be extrapolated from an input sample? Different solutions for the pdf estimates are produced when different log-likelihood target values are considered. In principle, target log-likelihoods should be randomized based on the pdf describing log-likelihoods for SURD (c.f. **Figure 2**). However, as a practical matter, the number of possible pdf estimates should be limited to keep computational cost down. As such, results from different fixed target values have been compared to determine if a single target value can be used.

A logical value for the target log-likelihood is the most probable value defined by the peak of the scoring function. This makes the computational aspect of the problem trivial. Among all possible pdf estimates based on different targets, the most probable log-likelihood target yields excellent visual results. Note that approximately 70% of all log-likelihoods are less than the most probable log-likelihood. Another reasonable target corresponds to 50% coverage. After much experimentation, and taking into consideration computational cost tradeoffs, it was found that a target log-likelihood value of $-0.37$ with approximately 40% coverage gives very good visual results for all test cases considered in this report. By dropping from 70% to 40% coverage, computing time was improved by worthwhile factors ranging from 10% to 50% depending on feature sharpness in the true pdf. It is interesting that the target value of $-0.37$ is slightly below the sample size invariant mean. Dropping below 40% coverage did not yield substantive benefits in computing time. Moreover, *similar* visual results are obtained when the coverage from a target value falls in the range from 20% to 60%. As such, 40% is used as the default coverage, and although an arbitrary number, it falls within a non-sensitive zone that consistently produces good quality solutions in an efficient manner for all cases considered in this report.

## 4. COMPUTATIONAL METHOD

We being this section by summarizing several components that have been described thus far: Given a random sample, $\{V_k\}$, the computational method performs two isomorphisms. The first isomorphism is static, in that the input data is rescaled by a one-time linear transformation to map $\{V_k\} \rightarrow \{X_k\}$ such that $-1 \leq X_k \leq 1$. The second isomorphism is dynamic to solve an inverse problem. A point $\lambda$ in parameter space of arbitrarily dimension controls the shape of the pdf estimate, $p_e(x|\lambda)$. For each trial pdf, an adaptive integration is carried out to calculate its associated cdf given by $F_e(x|\lambda)$. The trial cdf is then used to transform $\{X_k\} \to \{V_k\}$ where $U_k = F_e(X_k|\lambda)$. If the correct pdf is found, the statistical properties of the single order statistics are given by $\{U^{(s)}\}$ and they will be consistent with that for SURD. A universal scoring function based on a log-likelihood metric is used to test if the proposed single order statistics resembles what is typically expected from SURD.

The universal scoring function that is sample size invariant, shown in **Figure 2**, measures the log-likelihood of how far the single order statistics, $\{U^{(s)}\}$, deviates from SURD. The parameter space is explored by a random search method called funnel diffusion [33] to improve the transform $\{X_k\} \to \{U_k\}$ based on reaching 40% coverage in log-likelihood. The objective is to obtain a score that is near 40%. As such, statistical inference is extracted from the sampled data without over or under fitting because of the model independent nature of the scoring function and stopping criteria, both being invariant with respect to sample size. The implemented method produces an ensemble of solutions, each consistent with typical fluctuations from a finite size sample. A particular model pdf is a member of an ensemble of possible solutions. While any number of potential solutions can be generated to estimate uncertainty, the default setting used here for all cases is 5. The solution that is most similar to all other solutions is selected as the pdf estimate that is shown as the final result. In the following subsections, key parts of the computational method are briefly explained. A schematic flow chart for the entire computational method is shown in **Figure 3**.

**4.1 Iteratively Increasing Likelihood**

Any method that optimizes parameters of a function should suffice in increasing the log-likelihood score. Funnel diffusion is the random search method employed here. Funnel diffusion is conceptually described in this subsection with pseudo code given in Appendix A.1 for completeness. Funnel diffusion simulates a diffusion process implemented as a random walk in the parameter space, where the initial step size probes large length scales followed by a slowly decreasing step size to probe ever-smaller length scales as the solution point in the parameter space is approached. A random step is drawn from a Gaussian distribution with standard deviation, $\sigma_{FD}$, which gives the average step size. A random step corresponds to a trial move. If the log-likelihood given by Eq.(9) decreases, the trial move fails, and the walker does not move from the last position. The number of failures is counted as $n_F$. If the number of failures exceeds a maximum allowed number, $F_m$, the step size of the random walker is reduced such that $\sigma_{FD} \to f\sigma_{FD}$ where the multiplicative fraction $f$ sets an annealing rate. Each time the step size is reduced, the count for successive failures is reset to zero, and the random walk continues on a smaller length scale. Funnel diffusion stops when error in the parameters is sufficiently small, which is controlled by setting a minimum sigma. An advantage of using a random search method is that while all model pdf solutions generated in the ensemble are consistent with the MEM form, the random element of plausible solutions in the ensemble reflects typical sample size fluctuations.

A complication not encountered in random search methods appears in this application. The number of independent parameters (Lagrange multipliers) defines the dimension, $D$, of the search space, which is unknown. Recall parameter $\Lambda$ is a dependent variable fixed by the normalization condition on $[-1,1]$, determined after each random step. If funnel diffusion finishes with no solution found, $D$ is increased. As implemented here, $D = 0$ at first, then $D$ is incremented by 1 until it reaches 5, and thereafter increased by increments of 2. That is, $D$ is incremented as $0, 1, 2, 3, 4, 5, 7, 9, 11$ and so forth. The current parameter space defines a subspace once $D$ is increased, and the random walker is no longer confined to that subspace. Over many funnel diffusion cycles the dimension of parameter space is allowed to grow, but not to decrease. When the true pdf has sharp features, $D$ will become large. To reach large values of $D$ rapidly, the increment in $D$ can be increased. Alternatively, the increment in $D$ can become larger as $D$ increases (accelerated titration of new parameters). Several variants of how to expand the dimension of the parameter space were explored. An increment-step size of 2 beyond $D = 5$ provides a simple rule that produced robust results on all cases studied. Improving the efficiency in search spaces of variable dimension is left as an open optimization problem.

There is heuristic logic for gradually adding dimensions when failed solutions are found. When $D$ is increased, the additional Chebyshev polynomials oscillate more rapidly because they have more zero crossing points on $[-1,1]$. More samples are needed to ascertain with certainty the number of zeros a function has on an interval. Perhaps as little as 10 random points per zero crossing is sufficient to justify employing rapidly oscillating Chebyshev polynomials. If the true pdf has a tail with no wiggles, due to orthogonality and completeness of the Chebyshev polynomials, the expansion coefficients will decrease at higher orders. As such, evaluation of expansion coefficients consecutively from small to large orders makes sense. Furthermore, if a certain dimension is insufficient to yield a solution, the expansion has not yet converged, implying more terms are needed. This procedure tends to yield solutions with the smallest number of parameters possible, and it is easy to implement.

As discussed in section 3.3, the target value of $-0.37$ for the log-likelihood corresponding to $40\%$ coverage is a good point to stop funnel diffusion. However, as implemented, funnel diffusion stops when the rate of improving the scoring function drops below a minimum rate. In other words, funnel diffusion has its own criteria to stop, related to how much progress it makes in the annealing process, or if it reaches the maximum level of precision desired. The funnel diffusion method fails when the log-likelihood found is below the target value. After the search space dimension is increased, it again becomes possible to find a log-likelihood that is greater than the target value. As such, a better log-likelihood than the target value is needed to stop the process completely. For example, $D = 3$ is needed for a random sample of data that is normally distributed. Upon first attempts with $D = 1$ and $D = 2$, the data will not be described well. However, with $D = 3$, funnel diffusion is likely to reach ~70% coverage. Consequently, a Gaussian process is extremely fast to estimate because a large gain in coverage occurs once the dimension of the space spans the necessary parameters. Conversely, small gains in coverage are likely when the number of dimensions to accurately represent sharp features in the true pdf is large. In these cases, the final coverage will typically be just above or at 40%. But, what happens if the target log-likelihood value is never reached as the dimension continues to be increased?

Certain conditions are defined to guarantee the search process will terminate, either successfully or not. The process is successfully stopped when the target log-likelihood value is reached. However, the process is stopped before the target log-likelihood is reached when funnel diffusion makes too slow progress (stalls out) or if the parameter space exceeds a maximum dimension. When premature stopping occurs, a secondary condition declares the search successful if the log-likelihood is above $-1.166$, corresponding to 5% coverage. Otherwise, an error message declares no solution is found. Here, 5 different solutions are sought, but they must be accomplished successfully within 12 attempts before the program terminates. In all cases reported here, 5 solutions were successfully obtained in 5 attempts.

**4.2 Adaptive Integration**

Calculating the log-likelihood score using Eq. (9) and the penalty term using Eq. (10) as a function of $\{U^{(s)}\}$ requires integrating a trial pdf to obtain its associated cdf. Adaptive integration is implemented to ensure high accuracy while keeping computational cost low. The method is conceptually described in this subsection with pseudo code given in Appendix A.2 for completeness. Suppose the true pdf is known, $U_k = F(X_k)$ where $F(x)$ is the exact cdf, and let the interval $[0,1]$ be divided up into $M$ equal bins. Since $\{U_k\}$ represents SURD, these points will spread uniformly. This implies that bins of equal size will on average contain $N/M$ samples. For similar relative accuracy between sample points across different sample sizes, $M$ should be proportional to $N$. For a small sample size, a minimum resolution denoted as $M_{min}$ is set. For the nominal number of bins we invoke $M = \min[\max(bN, M_{min}), M_{max}]$ as a heuristic formula where $b = 0.005$, with a minimum allowed value of $M_{min} = 200$ and a maximum allowed value of $M_{max} = 1{,}500$. In principle no maximum is necessary, but it prevents the calculations from becoming slow with little loss of accuracy. For all test-distributions considered here, 1,500 served as a good nominal number of bins to carry out accurate numerical integration for up to $10^{20}$ data points. On $[0,1]$ the leading edge of each of the $M$ uniformly spaced bins is given by $u_k = \frac{k}{M}$, where $k$ runs from 0 to $M - 1$. Corresponding bin edges on the interval $[-1,1]$ are related by $u_k = F(x_k)$. This information is very important in setting up the adaptive integration, which requires an inverse mapping. The inverse mapping is called the quantile function, denoted by $Q(u)$. Then, $X_k = Q(u_k)$ and the $k$-th bin spans between $x_{k+1}$ and $x_k$. Note that bin size is not uniform for the $\{X_k\}$ variables, but each bin holds approximately $\frac{N}{M}$ number of data points.

In practice, non-uniform bins on $[-1,1]$ are determined by sorting $\{X_k\}$ values from smallest to largest. Based on $\lfloor N/M \rfloor$, a fixed number of samples are assigned to each bin. The sorted data is processed to locally partition data into bins by successively skipping over a fixed number of samples to delineate leading and trailing edges of bins. To maintain a minimum resolution, a maximum bin width is set to $2/M$ corresponding to a bin width of $1/M$ on $[0,1]$. If an identified bin is wider than the maximum width allowed, then it is subsequently partitioned into equal sized bins of smaller widths. Hence, $M$ is a nominal number of bins, and the interval $[-1,1]$ is divided into non-uniform bins that adapt to how data is distributed such that smaller bins are placed in regions of high accumulation of data, and wider bins are placed in data light regions. Simpson's rule for integration is applied to the non-uniform bins, offering high efficiency and ample accuracy.

### 4.3 Hierarchical Data Partitioning

Processing the sampled data hierarchically dramatically reduces computational cost. The strategy employed first sorts the full set of transformed data, $\{U^{(s)}\}$ on the range from $[0,1]$. The collection of these sorted data points has information about the global shape of the cdf that transforms $X_k$ to $U_k$. For example, out of 100,000 data points, the 27,000-th data point cannot be anywhere in the range from $[0,1]$, but rather, it must be after the 26,999-th point and before the 27001-th point. Therefore, uniformly selecting a subset of 1,000 points from 100,000 points by skipping every 100 points is not equivalent to throwing out 99% of the information. Rather, large-scale shape characteristics of the cdf are retained within this small subset of data. Sharper features are subsequently resolved as more uniformly spaced data is augmented. This physically intuitive idea is embodied in the expansion of orthogonal functions, starting with those that vary slowly to those with rapid variations (wiggles), as hierarchical partitioning successively augments more data points until all data points in the input file are considered.

The hierarchical process is invoked for datasets with more than 1,025 points; otherwise a single partition containing all the points is used. For larger datasets, the first step is to define a base partition having 1,025 points. Each partition thereafter has $(2^n + 1)$ points where, $n$ is incremented by 1, until the second largest partition is reached. The last partition contains all $\{U^{(s)}\}$ points. Each successively larger partition at the $n$-th level has one more point roughly halfway between two successive points from the previous partition at the $(n-1)$-th level. Note that these partitions are based on the sort order indexes. Only new points are added to larger partitions, such that a new partition contains all the points within previous partitions. When the log-likelihood score given by Eq. (9) was discussed in Section 3.2, the sample size being considered was the actual sample size, not the current partition size. To implement hierarchical partitioning, Eq. (9) has to be modified for each partition. Instead of using $N$ samples from the complete dataset, a smaller number for the partition called $N_p$ is used. The first term in Eq. (9) is found to be a self-averaging quantity for a given sample size provided the points are distributed uniformly, so the required modification is given as:

$$L\left(\{U^{(s)}\}_p\right) = \tfrac{1}{N_p}\ln\left(P_L\left(\{U^{(s)}\}_p\right)\right) - \tfrac{1}{2}\ln(N) \tag{11}$$

where $\{U^{(s)}\}_p$ is given by the subset of $N_p$ data points uniformly spaced across the sort ordered indices. The base partition has sufficient coverage of the sample to allow the search space dimension and the Lagrange multipliers to be rapidly determined. For a sample of 1 million data points, one calculation in the base partition will be almost 1000 times faster compared to using all the data. A single calculation in each partition thereafter takes about twice as long as it took in the immediate previous partition. It is worth mentioning that the sort ordering of the data is done only once, and all partitions observe the same original ordering of the sample dataset.

When moving into the next larger partition, the $\lambda$ parameters from the current partition are used as an initial condition. All partitions must be processed to complete the calculation. In cases where the true pdf does not have sharp features, all partitions beyond the first partition consume a tiny fraction of the total computational cost because the solution readily converges on the first partition. For example, convergence is fastest for a uniform pdf since it needs only one Lagrange multiplier. For cases where the true pdf has sharp features, the larger partitions are important to refine the estimates of Lagrange multipliers inferred from the lower level partitions. Sharp features require a larger number of Lagrange multipliers that slow down the calculations, but greater accuracy is gained as more rapidly oscillating orthogonal functions are incorporated in larger partitions. Although the hierarchical partitioning of data is not required to implement, it renders computational cost to an acceptable level for practical applications.

## 5. BENCHMARK EXAMPLES

The nonparametric MEM for pdf estimation is assessed on seven different types of test-distributions. Each type of distribution is presented in a different subsection to highlight how the approach handles interesting or difficult problems. In all cases, i.i.d. random samples were initially analyzed using *default settings* to perform the calculations with no a priori information about the data invoked. The pdf models obtained for each case are of good quality, although not all of these results are shown. To demonstrate versatility of the approach, other results are shown and discussed for illustrative cases where auxiliary constraints are imposed. For example, enforcing mirror-symmetry about $x = 0$ ensures a model pdf is an even function regardless of sample fluctuations, which increases model quality. Nevertheless, the dominant factor to achieve a good pdf estimate is the stopping criterion set by the threshold on the log-likelihood score (see section 3.3). In particular, SRQ-plots do not show significant improvement in their quality when auxiliary constraints are imposed, but computational time is often reduced.

**5.1 Generating and Assessing Results**

The same procedure is employed to conduct numerical experiments for all cases. Given a test-distribution defined by its pdf, $p(v)$, the corresponding quantile function, $Q(u)$, is constructed. A set of random samples $\{V_k\}$ are generated through the transformation $V_k = Q(r_k)$, where $r_k$ is a uniform random number on [0,1] using a high quality random number generator [62]. The scope of the testing involves four different sample sizes, with $N = 2^8, 2^{12}, 2^{16}, 2^{20}$, where the factors of 2 are considered only for convenience in presenting results. Performance does not depend on any special value of $N$. Due to greater statistical resolution inherent in larger sample sizes, the quality of the pdf estimate is monitored as $N$ increases. For each sample size, four independently drawn samples are randomly generated from the true pdf to produce a total of sixteen sample sets. Five possible model pdf solutions are calculated per input sample, yielding 20 model pdf solutions per sample size, and 80 model pdf solutions analyzed for each of the eight test-distributions.

For a given data sample, the central model pdf among 5 possible solutions is identified as the *pdf estimate*, having a minimum total pairwise squares error between it and the other four model solutions. The standard deviation for the differences between the alternate pdf models with respect to the central model (the estimate) can be used to report error bars on the pdf estimate (not shown here). An interesting question is how similar are the various pdf estimates (one for each sample) to each other? While pdf estimates for different samples can differ from the true pdf due to finite size fluctuations inherent within a sample, the method attempts to resist over fitting to sampling fluctuations. To see how well the method works for each test-distribution, the pdf-estimate is shown for each of the four input samples for each of the different sample sizes, and they are compared to the true pdf used to generate the samples. In addition, for each pdf estimate its corresponding SRQ-plot is shown.

For each test-distribution, four statistical metrics are tabulated to summarize the quality of the pdf estimate for each independent sample drawn from the true pdf. The metrics include the $p$-value for the Kolmogorov-Smirnov (KS) test [63] that compares the pdf estimate to the true pdf to determine if they are different or not. A $p$-value that is lower than 0.05 suggests the pdf estimate is different from the true pdf, although 5% of the time samples from the true pdf will return a p-value of 0.05 or less. The second metric is the Kullback–Leibler distance (KL-distance) [64, 65]. A *Figure of Merit* ($FOM$), defined in this work for the first time, provides a third metric. The $FOM$ is calculated by first generating an ensemble of samples from the estimated pdf that are the same size as the input sampled data. A series of pairwise comparisons are made between the generated samples to the input sample to quantify how likely the estimated pdf can produce the input sample data. The $FOM$ returns a number on the interval $(-\infty, 1]$, where a positive value indicates that the input sample has similar characteristics to a typical member in the ensemble. When $FOM > 1/2$ there is very high confidence that the input sample reflects a typical outcome from the estimated pdf. Details are given in Appendix A.3 for how the $FOM$ is calculated. The SURD level of coverage is reported as the fourth metric. In addition, the number of Lagrange multipliers used to determine the pdf is reported, which highlights variations between pdf estimates.

As discussed in section 3.3, the target coverage level for SURD is set to $40\%$. The SURD coverage metric is used to arrive at an estimated pdf, while the other three metrics defined above are used for independent benchmarking. Keeping in mind that $40\%$ is a target, it can happen that the SURD target is not reached after the funnel diffusion random search method is terminated, or the target is exceeded before the stopping criterion is checked. Based on SURD coverage *alone,* all pdf estimates reported here are of good quality. **Figure 4** shows the relationship between the KS $p$-value and the $FOM$ metric across the aggregate data from all test-distributions considered. This scatter plot shows that about $15\%$ of the pdf estimates reported here have either a low KS $p$-value and/or low $FOM$ value.

From **Figure 4** it is seen that when $FOM > \frac{1}{2}$ the corresponding KS $p$-values are typically found to be greater than ¼. Combined with good SURD coverage, consistent agreement among these three metrics provides a high level of confidence that the input sample is likely drawn from the estimated pdf. The KS $p$-value, $FOM$, and SURD metrics deal with probabilities for samples drawn from the same population distribution, and therefore significant differences *can occur* even if the estimated pdf is identical to the true pdf. This highlights the advantage of using multiple metrics to evaluate a model pdf. Furthermore, a consensus is not needed because outliers can occur without being indicative of an error. Conversely, a consensus of poor metrics identifies a pdf estimate that is of poor quality. With the exception of the KL-distance, the metrics are calculated without knowing the true population distribution. For $85\%$ of the results presented here, these latter three metrics unanimously indicate a fair to good solution. Among a total of 128 cases considered, only 3 cases were found to have a poor KS $p$-value and poor $FOM$ value, which is a number that can be statistically expected. These results are respectable considering that the test-distributions pose diverse challenges, where, in two cases, the default KDE implementation in MATLAB failed, and in two other cases the KDE results were notably worse visually.

### 5.2 Uniform Distribution

Here, the true pdf is defined as: $p(v) = \frac{1}{2}$ on $[-1,1]$. The uniform distribution serves as the simplest example of how the nonparametric MEM performs and, interestingly, it presents notable difficulties for some KDE methods due to its sharp boundaries, even when the boundaries are known a priori. For the results shown here, boundary values were not specified. A critical component of the MEM approach developed here is its ability to map sample data using a trial cdf, and *recognize* when the output from a trial cdf is equivalent to SURD. Since SURD is used to generate the uniformly random test samples, in principle the MEM should immediately find an exact solution using a single Lagrange multiplier (i.e. $\lambda_0$). In practice, this depends on the SURD target that is selected. For a low SURD target of 5%, it was found that indeed all pdf estimates were based on a single Lagrange multiplier, and thus each estimate was exact. Interestingly, by demanding greater accuracy using a larger SURD target percentage, such as 40%, additional Lagrange multipliers are typically invoked. However, it is the specific nature of the uniform distribution that afforded greater accuracy when the threshold criterion is relaxed. In general, as observed in all other cases, a 40% SURD target percentage is much better than a 5% target.

**Table 1** shows that the pdf estimate for all four samples at each size have ample quality. **Figure 5** shows the results visually, and it is seen that as the sample size increases the pdf estimates become flatter (closer to being uniform). In comparison to **Figure 1**, the SQR-plots of **Figure 5** show deviations within the same order of magnitude, and this visualization is an effective way to show that the estimated pdf is of high quality, having expected typical fluctuations. It is worth mentioning that the appearance of a higher density of data points plotted in **Figure 5**, compared to **Figure 1** for large sample sizes, is caused by the cap in the nominal number of integration bins given by $M_m = 1500$. Increasing $M_m$ can make SQR-plots look more similar to those in **Figure 1**, but the added time in the calculations (as much as a factor of 4 in some cases) did not warrant increasing the 1500 heuristic default value.

### 5.3 Laplace Distribution

Here, the true pdf is defined as: $p(v) = \frac{1}{2} e^{-|v|}$ on $(-\infty, \infty)$. This Laplace distribution provides an example of a function that cannot be expressed exactly as an exponential of a series expansion as in Eq. (2). Specifically, the cusp-discontinuity at $v = 0$ becomes rounded, and in practice this systematic error is within statistical resolution. However, there could be a physical reason why the pdf is known to be an even function, yet its shape not known. In this case the cusp can be modeled exactly by using symmetry as an auxiliary global constraint, where the data for $v < 0$ is imaged about the symmetry line at $v = 0$ so that all the data falls on $[0, \infty)$. After this transformation is made, the cusp that appears in this problem can be mathematically described exactly because the final answer is obtained by unfolding the estimated pdf back across the symmetry line for $v < 0$. Then the estimated pdf will be a continuous function and can have a discontinuous first derivative just like the true pdf.

An exact solution should contain two Lagrange multipliers $\lambda_0$ and $\lambda_1$. However, as **Table 2** shows, many more Lagrange multipliers appear in some solutions. For example, two pdf estimates invoke 19 different Lagrange multipliers. Not counting $\lambda_0$ corresponding to an overall normalization factor, from $\lambda_1$ to $\lambda_{19}$ the ratio $|\lambda_1|/\sum_{j=1}^{19}|\lambda_j|$ is greater than 99%. This means the high order Chebyshev polynomials in aggregate contribute less than 1% to the estimated pdf. While various pdf models deviate from the true pdf, they represent viable solutions and can be used to generate error bars for the pdf estimate. **Figure 6** shows that the pdf estimates are of high quality at all sizes, and the SQR-plots show the estimates are consistent with expected typical fluctuations. It is worth noting that better solutions could be obtained by running the funnel diffusion random search method longer for each dimension (i.e. number of Lagrange multipliers) before adding another dimension. This additional computational cost decreases the number of extra Lagrange multipliers and improves the quality of solutions. Nevertheless, all funnel diffusion parameters are fixed in this work to values that yield a good balance between accuracy and speed.

**5.4 Gamma Distribution**

Here, the true pdf is defined as: $p(v) = \frac{1}{\sqrt{\pi v}} e^{-v}$ on $(0, \infty)$. This special case of a gamma distribution contains a square-root singularity, and it mimics van Hove singularities that appear in the density of states of free electron models in condensed matter physics. Important to this work is to test the method on a pdf with a singularity. Although modeling a singularity is challenging for KDE, the MEM approach described here requires no special treatment. Despite the divergence when $v \to 0$, the functional form of Eq. (2) provides a good representation of the sampled data within statistical resolution. A larger number of Lagrange multipliers are needed to model the singularity with accuracy commensurate to statistical resolution as sample size increases. This effect is clear in **Table 3**, where it is seen that as sample size increases the divergence is modeled more faithfully by using more Lagrange multipliers. The pdf estimates are shown in **Figure 7** where no boundary conditions are specified.

The $FOM$ metrics and the visual appearance of systematic error in the SQR-plots together indicate the pdf estimates have some flaws. In particular, it is seen that for small sample size the pdf estimates have smooth large-scale wiggles that are not from random noise. When the number of Lagrange multipliers is small the pdf estimates incur systematic error because the true pdf cannot be represented by a truncated series of Chebyshev polynomials that has not yet converged. However, as sample size increases more Chebyshev polynomials are used, and this relegates the systematic wiggles to appear on a finer length scale. For the largest sample size, systematic error is not visible in the pdf estimate, but the SQR-plot shows deviations with overall smaller amplitude than the expected typical fluctuations. Moreover, systematic errors show up as a regular pattern near the singularity, and as spikes in the tail. This means random noise is not displayed in the SQR-plot as it should be. Of course it is important to realize that within the context of the critical analysis presented here, a traditional QQ-plot would suggest the estimated pdf is perfect, and a KS $p$-value test flares well. Test cases such as this one motivated the development of the $FOM$ metric, which is very sensitive to local deviations in probability density. This test case highlights a general result that errors due to random fluctuations from finite size samples and systematic errors due to inadequate functional representations compete with one another during the minimization of the log-likelihood scoring function.

Regardless of the basis functions used, when a truncation is employed in a series expansion, error in representing the true pdf will be incurred that is not related to statistical resolution. Adding additional Lagrange multipliers is the only way to reduce this systematic error. For example, it is possible for a user to enforce using a sufficient number of Lagrange multipliers to produce model solutions with virtually no systematic error. However, in doing this, less random error also occurs, and the data is over fitted with SURD coverage approaching 100%. Therefore, the approach taken here of adding a minimum number of Lagrange multipliers as SURD coverage increases affords a way to be approximately consistent with expected levels of fluctuations, but these fluctuations reflect systematic and random errors combined. The artifact of systematic error remains an open problem within the current implementation. Fortunately, the SQR-plots identify these systematic errors and highlight the overall quality of an estimated pdf. Reducing systematic errors may require including additional SURD characteristics as part of a target, such as binning statistics. Nevertheless, it is demonstrated here that the employed implementation provides a robust automated process that yields good results when considering the balance between automation, speed and accuracy without over-fitting the data.

**5.5 Sum of Two Gaussian Distributions**

Here, the true pdf is defined as: $p(v) = \frac{7}{10}\mathcal{N}(v \mid \mu_1 = 5, \sigma_1 = 3) + \frac{3}{10}\mathcal{N}(v \mid \mu_2 = 0, \sigma_2 = \frac{1}{2})$ on $(-\infty, \infty)$ where $\mathcal{N}(v \mid \mu, \sigma)$ denotes a Gaussian distribution with mean $\mu$ and standard deviation $\sigma$. This binary mixture of two Gaussian functions produces random samples with bimodal character, and highlights a case where two markedly separated spatial resolution scales are present. This example was contrived so that the pdf cannot be simply expressed as a single exponential in accordance with Eq. (2). The less probable population appears as a small shoulder protruding out from the side of a broad Gaussian form. A large number of Lagrange multipliers will be necessary to match the functional form of the true pdf. By design, small sample sizes do not have the statistical resolution necessary to discern the difference between a sample fluctuation and the presence of a sharp Gaussian shoulder hiding within the broader Gaussian. The results in **Table 4** show that the number of Lagrange multipliers increase as sample size is increased, and the quality of the pdf estimate correspondingly improves as shown in **Figure 8**. Note that the small shoulder is essentially ignored at a sample of size of 256 data points because this MEM method resists over-fitting data. Similar to the gamma distribution results, **Figure 8** shows systematic error in the form of oscillations in the SQR-plot for all but the smallest sample sizes. Although this case is easily handled by KDE methods, the small protruding peak is over fit at small sample size.

**5.6 Five Fingers Distribution**

Here, the true pdf is defined as: $p(v) = w \sum_{k=1}^{5} \frac{1}{5} \mathcal{N}\left(v \mid \mu_k = \frac{2k-1}{10}, \sigma = \frac{1}{100}\right) + (1-w)$ on $[0,1]$ where $\mathcal{N}(v \mid \mu, \sigma)$ denotes a Gaussian distribution with mean $\mu$ and standard deviation $\sigma$. It is seen that 5 sharp Gaussian distributions are added to a uniform distribution. Conceptually, the probability, $w$, is taken from the uniform distribution and subsequently redistributed equally to five Gaussian distributions. Two cases are considered here, corresponding to $w = 0.5$ and $w = 0.2$. Due to the sharpness of all five Gaussian distributions it is virtually exact to consider the domain to be on the range $[0,1]$ without cutting off the tails of the Gaussian distributions and without having to renormalize the pdf. The five fingers distribution demonstrates how the nonparametric MEM resolves the Gaussian shaped fingers better as statistical resolution is increased with greater sample size. Moreover, the larger fingers ($w = 0.5$) are resolved more easily than the smaller fingers ($w = 0.2$).

As seen in **Table 5** and **Table 6** for large and small fingers respectively, more Lagrange multipliers are needed to resolve features that are commensurate with sample size. For the same sample size, the smaller fingers require a smaller number of Lagrange multipliers indicating they are resolved more easily than the larger fingers. **Figure 9** and **Figure 10** show that the large and small five-finger-distributions are well described. In contrast, poor results were obtained for the five-finger-distributions using KDE, especially for large sample sizes. Importantly, at a sample size of 256 points the nonparametric MEM cannot discern the fingers compared to random sample fluctuations, which means the estimated pdf is close to a uniform distribution (more so for $w = 0.2$). The oscillations seen in the SQR-plots in **Figures 9** and **10** exemplify the compromise made between systematic and random errors as discussed above.

Indeed, the pdf estimates have the appearance that is reminiscent of using a truncated Fourier series to approximate a periodic waveform with sharp features.

## 5.7 Cauchy Distribution

Here, the true pdf is defined as: $p(v) = \frac{b}{\pi(v^2+b^2)}$ on $(-\infty, \infty)$ where $b = \frac{1}{2}$. This $p(v)$ is the classic Cauchy distribution, also called the Lorentzian distribution, which has a heavy tail and models extreme statistics in the sense that the second moment, $\langle v^2 \rangle$, and the standard deviation do not exist. For a finite sample, one can calculate the variance and standard deviation, but as more samples are drawn from this distribution the standard deviation and variance continue to grow indefinitely. Distributions with extreme statistics are generally challenging for a KDE approach, unless the tails are severely cut in advance, but this creates error in the estimate. It turns out that the MEM presented here requires censor windowing for it to provide robust and stable results. The employed default censor window on $[a, b]$ are assigned limits given by $a = Q_{25} - 7(Q_{75} - Q_{25})$ and $b = Q_{75} + 7(Q_{75} - Q_{25})$ and this range is sufficiently restrictive to handle this Cauchy distribution without user-defined overrides. All calculations are stable and the results within the censor window do not depend on the data that is filtered out. As such, this censor windowing boundary rule in terms of quantiles is applied to any arbitrary input data sample for which the user has no knowledge.

In **Figure 11**, it is clear that increasing the number of samples provides a better representative of the Cauchy distribution, and the SQR-plots show less systematic errors than the other difficult distributions considered here. Despite some wiggles that appear in the tails of the pdf estimates at small sample sizes, the nonparametric MEM provides robust statistical inference. For the $N = 256$ sample size It is seen in **Figure 11** that the tails in the pdf estimates have bumps that decay to zero further out. Although the bumps in the tails at small sample size certainly do not look good, it is consistent with SURD at the set target value. It is likely the bumps in the tails are reflecting what would surely be considered outliers if the data was viewed as coming from a Gaussian distribution. If some smoothing is applied, it is often possible to obtain a model pdf that looks identical to the true pdf for some cases, but inevitably other cases become worse. The best compromise we found is not to smooth at all, even though smoothing remains a user-defined option following past work [33]. With user-defined options, an experienced analyst can perform a series of trial and error adjustments to prevent over or under smoothing by eye. However, subjective decisions by a user of statistical software are what this nonparametric MEM avoids. Therefore, the nonparametric MEM presented here makes no attempt at performing smoothing.

## 5.8 Discontinuous Distribution

Here, the true pdf is defined as a discontinuous distribution on the interval [0,1], given as:

$$p(v) = \begin{cases} 4/5, & v < 0.3 \text{ or } v > 0.8 \\ 1, & v > 0.4 \text{ and } v < 0.5 \\ 5/4, & \text{otherwise} \end{cases}$$

The expansion in Eq. (2) can approximate a discontinuous distribution within the accuracy of statistical resolution. For small sample sizes, infinitely sharp discontinuous features are replaced by broad smooth differentiable continuous curves that make no attempt to follow the abrupt changes. As sample size increases the discontinuous features are automatically modeled more accurately as statistical resolution increases. For sharp edges to come into focus the number of Lagrange multipliers must be increased. **Table 8** summarizes the characteristics of the pdf estimate and **Figure 12** shows the pdf estimates and the SQR-plots. Comparison of the pdf estimates with the true pdf shows that the sharp discontinuities are captured, albeit more wiggles are found along the flat plateaus than what KDE produces. The large systematic errors observed in the SQR-plots reflect rounding of the discontinuities.

For small samples, results from this nonparametric MEM are very similar to those obtained by a KDE method. At large sample sizes, KDE methods provide smoother results than the MEM used here due to systematic error as discussed above. It is worth mentioning that for this distribution (and other distributions that yield large systematic errors) the pattern that appears in a SQR-plot is virtually the same across all independent samples of a given size. This is a strong signature of systematic error. Unless it is possible to include more Lagrange multipliers without over-fitting the data within statistical resolution, the user can use the SQR-plots to elucidate were errors occur. For example, by careful inspection of **Figure 12**, it can be seen that positive peaks and negative troughs in SQR-plots respectively reflect over and under estimates of probability density relative to the true pdf. Because of this direct relationship, residual plots can be leveraged for further analysis. Given that systematic error in the pdf estimates is distinct from random noise (c.f. **Figure 1**), a fruitful direction we have explored is to employ signal processing methods to extract systematic deviations from a quantile residual plot and use them to correct the pdf estimate. Although this approach does improve the pdf estimates, the gains are related to smoothing out fine scale wiggles and do not change the overall performance characteristics of the non-parametric MEM. Although beyond the scope of this report, post signal processing SQR-plots shows promise as a means to obtain better pdf estimates without subjective human intervention. This additional signal processing feature will be described elsewhere.

**5.9 Run Times as a Function of Sample Size**

The average computational times to calculate a model pdf for the distributions that were considered in this report are benchmarked in **Figure 13** for different sample sizes ranging from 256 to 1,048,576 data points. The average computational times were obtained from calculating 100 independent solutions for each sample size, 25 solutions for each of the four different data sets. The computational cost empirically scales nearly linearly over three orders of magnitude in the number of sample data points. The slowest case among all the distributions discussed is the gamma distribution, taking on average 26 minutes on a modern laptop to obtain a model pdf for $2^{20}$ data points. This time is necessary to model the singularity correctly. The next slowest case is the Cauchy distribution, which takes an average of 2.5 minutes for $2^{20}$ data points. Recall that both of these distributions have a rapid increase in number of Lagrange multipliers as the number of sample points increases.

Hierarchical partitioning of the sampled data makes the cost of calculations tractable, and faster than KDE in some cases. Based on the structure of the algorithm (c.f. **Figure 3**) and benchmarked timings, a model pdf can be calculated for a sample containing $2^{20}$ data points in well under 1 minute on a modern laptop for distributions requiring few Lagrange multipliers, such as uniform, Laplace and Gaussian. For any distribution, sample sizes of 25,000 data points will on average take less than 1 minute to calculate a model pdf. This is because for difficult distributions that need a large number of Lagrange multipliers at high statistical resolution, a sample size of 25,000 data points is too small to justify using a large number of Lagrange multipliers. For sample sizes of less than 256 data points, poor statistical resolution only allows pdf estimates to be very smooth, which keeps computational cost minimal with net speed that is competitive to KDE methods. Compared to KDE on difficult cases where KDE succeeded, the presented MEM can be more than an order of magnitude slower. It is worth noting that the nonparametric MEM presented here has not been optimized to the same level as highly efficient KDE methods with decades of development. With room for further optimization the nonparametric MEM introduced here is likely to become more competitive in computing times while solving troubling areas that inflict KDE methods.

**6. CONCLUSIONS**

We developed a nonparametric maximum entropy method to estimate a pdf for a sample of random data where virtually no information about the data is known. Of course, any information about a random process or constraints on the data should be leveraged to obtain the best possible model. In this paper, we assume observations are independently and identically distributed and that the pdf is a continuous and infinitely differentiable function that yields maximum entropy consistent with the observations that impose constraints on its shape. The maximum entropy assumption dictates the functional form of the pdf. A random search method is employed to determine the coefficients of a series expansion over a set of orthogonal basis functions. A trial cdf transforms sampled random data, $\{V_k\}$, onto the interval $[0,1]$ via $U_k = CDF(V_k)$. The new random variables, $\{U_k\}$, are subsequently assessed using a scoring function based on log-likelihoods of single order statistics of *sampled uniform random data* (SURD). The scoring function is universal, being problem independent and sample size invariant. Plausible cdf solutions are identified when the proposed $\{U_k\}$ have typical SURD characteristics. As an inverse problem, there is no unique model pdf for a finite size sample. Therefore, an ensemble of model pdf solutions is generated from which a typical model pdf is selected as a pdf estimate for the population. The variation among plausible model pdf solutions from this ensemble is used to place uncertainties on the pdf estimate. The analytical form of the pdf estimate approximates the true pdf markedly well within statistical resolution. A major advantage of this approach is that it is resistant to over and under fitting to sample fluctuations.

Many technical difficulties were overcome to arrive at a numerically stable and tractable method. For example, a process was developed in which the number of parameters to optimize is not predetermined. Rather, the dimension of parameter space expands systematically. Parameters correspond to Lagrange multipliers, and they serve as coefficients in a series expansion involving an orthogonal basis set of functions. Titrating Lagrange multipliers into the series expansion to describe increasingly fine features makes sense intuitively because the shape of a model pdf is resolved on larger length scales before smaller length scales. Consequently, more data points are needed to justify extending the series to more terms in order to capture sharper features on small length scales. Furthermore, augmenting data hierarchically as smaller scale features are refined increases the speed of the calculations dramatically. This method has been demonstrated to work well on a challenging test set of distributions, and it is also shown to be a versatile approach for specifying boundary conditions, symmetries and censor windows.

The key notion leveraged in this approach is characterizing typical fluctuations. Typical fluctuations scale with sample size and have certain characteristics that are described by the universal scoring function for SURD. Deviations from typical SURD fluctuations are easy to see using scaled quantile residual plots. It is seen that when the series is truncated prematurely, meaning before it converges, systematic error competes with random error. The scaled quantile residual plots are sensitive to errors at the resolution of statistical noise because they are scaled by a factor of $\sqrt{N+2}$ compared to a standard residual plot. As such, the current method provides a reasonable estimate for a pdf that respects the interplay between statistical inference and statistical resolution. Nevertheless, a weakness of this approach is that, as the log-likelihood function is minimized, systematic error is reduced at a slower rate than random error as the total error becomes consistent with the typical amplitude that is expected in SURD fluctuations. Many avenues are available to address this weakness and to improve the method. These include signal processing to extract signal from noise, monitoring convergence rates to achieve better accuracy, including binning statistics as part of the scoring function criteria, and optimizing the funnel diffusion random search method with more responsive stopping criteria to achieve better efficiency. As we pursue these directions the Java program used to generate the data in this report (without modification) can be obtained by request. The initial release of the Java program under a general public license is planned, and thereafter, we plan to provide periodic improvements.

In the modern era of data analytics, error free data driven automation for high throughput analysis applications is essential. The nonparametric MEM introduced here is shown to exhibit robust solutions that are consistent with one another. The method itself is very versatile, as demonstrated by a range of different numerical experiments. Without invoking any subjective human judgments that are required for other methods, an automated nonparametric estimate for a pdf in analytic form is obtained without over or under fitting to the data, and uncertainties are quantified. The nonparametric MEM presented here provides a holistic approach to statistical inference. To quantify the quality of a model pdf, a large part of our results rely on the universal scoring function. We made use of the scaled quantile residual plots and introduced the figure of merit for independent assessments. These statistical metrics are applicable to other estimation methods, including KDE. In particular, the scoring function can be used to set objective criteria to resist under and over fitting to a random sample of data.

**ACKNOWLEDGMENTS**

We thank Dr. Michael Grabchak for many fruitful discussions on our approach and applications to extreme statistics, and thank the Center for Biomedical Engineering and Science for financial support given to DJJ.

**APPENDIX – Pseudo Code**

**A.1 Funnel Diffusion**

Listed below are the detailed steps for the main program logic, which explains how funnel diffusion works together with a trial cdf and scoring method to converge toward a final solution. Each of the three main sections (initialization, main loop, and end processing) are executed a user-defined number of times (default is 5), in order to produce an ensemble of viable model pdf solutions.

**Initialization Steps**

1. Initialize one Lagrange multiplier with a value of zero
2. Create pdf for a uniform distribution, calculate and normalize the cdf
3. Map transformed sample data to [0, 1] using cdf
4. Using order statistics, create a log-likelihood score for mapped data
5. If log-likelihood score is consistent with SURD within user-defined target, end program and accept uniform distribution as a solution (default target is 40%)
6. Add an additional Lagrange multiplier and initialize it to zero
7. Set current funnel diffusion step size, sigma, to user-defined value (default is 0.1)
8. Set the funnel diffusion decay_rate to user-defined value (default is $\frac{\sqrt{2}}{2}$), such that sigma → decay_rate*sigma.
9. Set final allowed step size to user-defined value (default is 0.001).
10. Execute Main Loop below until one of the following three exit criteria occurs:
    a. Acceptable solution is found
    b. User-defined maximum number of Lagrange multipliers is reached (default max = 300)
    c. Score has not improved by a minimum (user-defined) percentage in a user-defined number of consecutive Lagrange multiplier additions (defaults are 1% and 3, respectively)
11. Initialize funnel-step loop counter to zero

**Main Loop**

1. Execute funnel diffusion on all Lagrange multipliers by making a random step using a Gaussian distribution from each of the current values, using the current sigma as the standard deviation
2. Create a pdf using the new Lagrange multipliers, calculate and normalize the cdf
3. Map transformed sample data to [0, 1] using cdf
4. Using order statistics, create a log-likelihood score for mapped data
5. If log-likelihood score is consistent with SURD within defined target, end program and accept current pdf as the solution, else continue to step 6

6. If funnel-step loop counter < user-defined maximum, increment loop counter by one, and go to step 1 (default for max loop counter is $F_m = 100$), else continue to step 7
7. If current sigma is > final allowed step size, decrease current sigma by decay_rate, reset loop counter to zero, and go to step 1, else continue to step 8
8. If number of Lagrange multipliers < maximum number allowed, add an additional two Lagrange multipliers, initialize them to zero, increment loop counter by one, and go to step 1, else continue to step 9
9. Reassess current log-likelihood score to determine if it meets a minimum SURD threshold of 5%. If the score meets this lower threshold, tag this solution as a success, else tag solution as a failure
10. End program

**End Processing**

1. Apply the Lagrange multipliers on raw input data to create pdf and normalized cdf for final solution.
2. Write out full solution to a file, even if it was a failure, including all parameters, score, and raw data with pdf and cdf

**A.2 Equal Partitioning of Sample Quantiles**

The following steps outline the process of deriving a variable sized $dx$ value to use for calculating and integrating each trial pdf. The algorithm described below allows the Simpson's rule to be applied adaptively. Although the analytical function is not determined until the end, the data is used to determine the support for integration. Integration accuracy does not depend on the fluctuations that occur across data samples.

1. Transform raw data on interval [0, 1] and sort order
2. Calculate minimum number of integration points, denoted by $M$, as a function of sample size $N$: $M = 200 + N * 0.005$. If $M > 1500$, set $M = 1500$.
3. Calculate maximum allowed delta size; $\Delta x_{max} = \frac{2}{M-1}$.
4. Calculate bin size for number of sample points per integration point; $b = floor\left(\frac{N}{M-1}\right)$.
5. Set $\delta = data(b)$, and $k = b$, where $data(b)$ is the $b^{th}$ value of the sorted transformed data.
6. If $\delta x \leq \Delta x_{max}$, then add $\Delta x_{max}$ to the list of variable $dx$ values.
7. If $\delta x > \Delta x_{max}$ then add $floor\left(\frac{\delta x}{\Delta x_{max}}\right)$ equal values of $\Delta x_{max}$ to the list of variable $dx$ values.
8. Calculate the difference $\delta = data(k + b) - data(k)$. Set $k = k + b$.
9. Repeat steps 6, 7, and 8 for each $b$ group of sorted data points, in order from smallest to largest.

## A.3 Figure of Merit

The inputs to this algorithm include the original data sample, and random data samples that are generated from the cdf of the proposed solution. The steps taken to calculate FOM are as follows.

1. Create 10 data samples from the cdf, sort them and call these samples 'reference samples'
2. Create an additional 100 data samples from the cdf, sort them and call these samples 'test samples'
3. Calculate the average difference, per position k, between the original data sample and the 100 test samples. Call this $\mu_{sample}\{k\}$.
4. Calculate the average difference, per position k, between *each* of the 10 reference samples and the 100 test samples. Call this $\mu_{test}\{10,k\}$.
5. Find the average and standard deviation, per position *k*, of the differences between all reference samples and test samples, by comparing each 10x100 sample set. Call these values $\mu_{ref}\{k\}$ and $\sigma_{ref}\{k\}$.
6. For each of the 10 reference samples, calculate a FOM per position.

$$FOM_{\{10,ref\}} = 1 - \frac{|\mu_{ref}\{k\} - \mu_{test}\{10,k\}|}{2*\sigma_{ref}\{k\}}$$

7. Calculate a FOM per position for the original data sample.

$$FOM_{sample} = 1 - \frac{|\mu_{ref}\{k\} - \mu_{sample}\{k\}|}{2*\sigma_{ref}\{k\}}$$

8. Calculate the average and standard deviation, $\langle FOM_{ref} \rangle$ and $FOM_{\sigma_{ref}}$, averaged over the 10 reference samples and *k* positions.
9. Calculate the average, $\langle FOM_{sample} \rangle$, averaged over all *k* positions.
10. Using the values calculated in steps 8 and 9, compute a single number for FOM,

$$FOM = 1 - \frac{|\langle FOM_{ref} \rangle - \langle FOM_{sample} \rangle|}{2*FOM_{\sigma_{ref}}} .$$

**Tables with table captions**

| Sample Size | p-value | KL distance | Figure of Merit | SURD coverage | Lagrange Multipliers |
|---|---|---|---|---|---|
| 256 | 0.47 | 5.0E-6 | 0.89 | 0.61 | 1 |
| | 0.62 | 2.0E-4 | 0.88 | 0.42 | 2 |
| | 0.62 | 5.0E-6 | 0.76 | 0.90 | 1 |
| | 0.64 | 5.0E-6 | 0.85 | 0.61 | 1 |
| 4096 | 0.49 | 5.0E-5 | 0.90 | 0.51 | 1 |
| | 0.53 | 5.0E-6 | 0.83 | 0.62 | 1 |
| | 0.56 | 5.0E-6 | 0.85 | 0.53 | 1 |
| | 0.61 | 5.0E-6 | 0.85 | 0.85 | 1 |
| 65536 | 0.46 | 3.7E-5 | 0.91 | 0.45 | 3 |
| | 0.47 | 5.0E-6 | 0.71 | 0.54 | 1 |
| | 0.47 | 5.0E-6 | 0.80 | 0.56 | 1 |
| | 0.55 | 5.0E-6 | 0.89 | 0.76 | 1 |
| 1048576 | 0.52 | 6.7E-6 | 0.82 | 0.40 | 5 |
| | 0.53 | 6.3E-6 | 0.84 | 0.42 | 7 |
| | 0.56 | 6.4E-6 | 0.88 | 0.49 | 3 |
| | 0.59 | 5.0E-6 | 0.83 | 0.46 | 1 |

**Table 1:** Characteristics of pdf estimates for each of four independent samples drawn from a uniform distribution without boundary conditions specified.

| Sample Size | p-value | KL distance | Figure of Merit | SURD coverage | Lagrange Multipliers |
|---|---|---|---|---|---|
| 256 | 0.51 | 7.8E-3 | 0.87 | 0.47 | 2 |
|  | 0.58 | 1.3E-2 | 0.87 | 0.43 | 2 |
|  | 0.59 | 9.0E-3 | 0.87 | 0.40 | 4 |
|  | 0.68 | 1.0E-2 | 0.84 | 0.40 | 4 |
| 4096 | 0.45 | 4.8E-3 | 0.65 | 0.45 | 2 |
|  | 0.53 | 4.2E-3 | 0.85 | 0.45 | 2 |
|  | 0.54 | 5.3E-3 | 0.86 | 0.50 | 2 |
|  | 0.56 | 2.5E-3 | 0.85 | 0.46 | 2 |
| 65536 | 0.50 | 2.9E-3 | 0.88 | 0.40 | 9 |
|  | 0.50 | 4.3E-3 | 0.76 | 0.41 | 5 |
|  | 0.53 | 4.4E-3 | 0.82 | 0.43 | 15 |
|  | 0.62 | 4.4E-3 | 0.80 | 0.49 | 2 |
| 1048576 | 0.34 | 4.3E-3 | 0.58 | 0.40 | 9 |
|  | 0.57 | 4.3E-3 | 0.82 | 0.40 | 19 |
|  | 0.62 | 4.3E-3 | 0.83 | 0.40 | 2 |
|  | 0.68 | 4.3E-3 | 0.84 | 0.40 | 19 |

**Table 2:** Characteristics of pdf estimates for each of four independent samples drawn from a Laplace distribution. The symmetry condition that the pdf must be an even function is imposed. No boundary conditions were specified.

| Sample Size | p-value | KL distance | Figure of Merit | SURD coverage | Lagrange Multipliers |
|---|---|---|---|---|---|
| 256 | 0.59 | 9.2E-3 | 0.90 | 0.40 | 9 |
|  | 0.59 | 1.6E-2 | 0.89 | 0.40 | 9 |
|  | 0.60 | 1.6E-2 | 0.90 | 0.40 | 9 |
|  | 0.60 | 2.4E-2 | 0.88 | 0.41 | 13 |
| 4096 | 0.31 | 1.1E-2 | 0.59 | 0.40 | 33 |
|  | 0.33 | 1.0E-2 | 0.64 | 0.40 | 29 |
|  | 0.38 | 1.1E-2 | -5.98 | 0.40 | 31 |
|  | 0.41 | 1.1E-2 | 0.65 | 0.41 | 29 |
| 65536 | 0.17 | 1.0E-2 | -0.02 | 0.41 | 69 |
|  | 0.17 | 1.0E-2 | -0.23 | 0.43 | 77 |
|  | 0.17 | 9.9E-3 | 0.02 | 0.40 | 77 |
|  | 0.18 | 9.6E-3 | -0.15 | 0.40 | 81 |
| 1048576 | 0.28 | 9.9E-3 | -0.41 | 0.18 | 239 |
|  | 0.36 | 9.9E-3 | -0.37 | 0.17 | 233 |
|  | 0.44 | 9.9E-3 | 0.30 | 0.34 | 213 |
|  | 0.46 | 9.8E-3 | 0.15 | 0.28 | 289 |

**Table 3:** Characteristics of pdf estimates for each of four independent samples drawn from a gamma distribution that has a square root singularity.

| Sample Size | p-value | KL distance | Figure of Merit | SURD coverage | Lagrange Multipliers |
|---|---|---|---|---|---|
| 256 | 0.50 | 1.3E-2 | 0.90 | 0.43 | 3 |
| | 0.51 | 1.3E-2 | 0.89 | 0.46 | 3 |
| | 0.54 | 1.4E-2 | 0.90 | 0.40 | 4 |
| | 0.62 | 1.3E-2 | 0.87 | 0.42 | 3 |
| 4096 | 0.46 | 5.7E-3 | 0.92 | 0.41 | 17 |
| | 0.51 | 5.2E-3 | 0.92 | 0.41 | 17 |
| | 0.54 | 7.8E-3 | 0.91 | 0.42 | 13 |
| | 0.55 | 9.0E-3 | 0.91 | 0.40 | 13 |
| 65536 | 0.51 | 5.2E-3 | 0.89 | 0.40 | 33 |
| | 0.53 | 5.8E-3 | 0.87 | 0.40 | 33 |
| | 0.55 | 6.1E-3 | 0.89 | 0.40 | 33 |
| | 0.62 | 7.5E-3 | 0.84 | 0.42 | 29 |
| 1048576 | 0.43 | 6.7E-4 | 0.90 | 0.40 | 53 |
| | 0.50 | 7.2E-3 | 0.90 | 0.40 | 49 |
| | 0.51 | 5.1E-3 | 0.89 | 0.40 | 47 |
| | 0.57 | 5.0E-3 | 0.90 | 0.40 | 47 |

**Table 4:** Characteristics of pdf estimates for each of four independent samples drawn from a bimodal distribution described by a binary mixture model referred in this work as the two Gaussian distribution.

| Sample Size | p-value | KL distance | Figure of Merit | SURD coverage | Lagrange Multipliers |
|---|---|---|---|---|---|
| 256 | 0.33 | 3.4E-1 | 0.89 | 0.42 | 15 |
| | 0.37 | 4.1E-1 | 0.89 | 0.41 | 11 |
| | 0.39 | 4.7E-1 | 0.88 | 0.41 | 15 |
| | 0.42 | 4.1E-1 | 0.90 | 0.43 | 7 |
| 4096 | 0.17 | 1.1E-1 | 0.62 | 0.40 | 35 |
| | 0.20 | 9.5E-2 | 0.62 | 0.40 | 37 |
| | 0.20 | 1.1E-1 | 0.49 | 0.40 | 35 |
| | 0.22 | 1.3E-1 | 0.57 | 0.40 | 35 |
| 65536 | 0.06 | 1.6E-2 | 0.55 | 0.40 | 51 |
| | 0.07 | 1.5E-2 | 0.43 | 0.41 | 53 |
| | 0.07 | 1.5E-2 | 0.65 | 0.40 | 53 |
| | 0.08 | 1.6E-2 | 0.62 | 0.41 | 51 |
| 1048576 | 0.02 | 2.4E-3 | 0.04 | 0.42 | 99 |
| | 0.03 | 2.2E-3 | 0.32 | 0.42 | 101 |
| | 0.05 | 2.1E-3 | 0.87 | 0.42 | 97 |
| | 0.06 | 2.4E-3 | 0.84 | 0.42 | 101 |

**Table 5:** Characteristics of pdf estimates for a five fingers distribution with a 0.5 weight.

| Sample Size | p-value | KL distance | Figure of Merit | SURD coverage | Lagrange Multipliers |
|---|---|---|---|---|---|
| 256 | 0.58 | 8.8E-2 | 0.87 | 0.41 | 3 |
|  | 0.60 | 8.3E-2 | 0.93 | 0.42 | 3 |
|  | 0.61 | 8.6E-2 | 0.89 | 0.43 | 3 |
|  | 0.63 | 8.2E-2 | 0.91 | 0.40 | 3 |
| 4096 | 0.23 | 4.4E-2 | 0.82 | 0.41 | 15 |
|  | 0.25 | 4.2E-2 | 0.76 | 0.40 | 15 |
|  | 0.28 | 4.2E-2 | 0.76 | 0.40 | 15 |
|  | 0.33 | 4.2E-2 | 0.89 | 0.42 | 15 |
| 65536 | 0.13 | 1.1E-2 | 0.61 | 0.41 | 39 |
|  | 0.14 | 1.1E-2 | 0.51 | 0.40 | 39 |
|  | 0.17 | 1.1E-2 | 0.64 | 0.40 | 39 |
|  | 0.20 | 1.1E-2 | 0.55 | 0.40 | 39 |
| 1048576 | 0.10 | 1.5E-3 | 0.74 | 0.40 | 67 |
|  | 0.10 | 1.4E-3 | 0.71 | 0.40 | 67 |
|  | 0.11 | 1.4E-3 | 0.69 | 0.40 | 67 |
|  | 0.12 | 1.4E-3 | 0.78 | 0.40 | 69 |

**Table 6:** Characteristics of pdf estimates for a five fingers distribution with a 0.2 weight.

| Sample Size | p-value | KL distance | Figure of Merit | SURD coverage | Lagrange Multipliers |
|---|---|---|---|---|---|
| 256 | 0.49 | 5.1E-3 | 0.87 | 0.42 | 9 |
| | 0.53 | 3.8E-3 | 0.85 | 0.41 | 9 |
| | 0.54 | 4.0E-3 | 0.82 | 0.41 | 9 |
| | 0.63 | 4.6E-3 | 0.86 | 0.41 | 9 |
| 4096 | 0.32 | 6.5E-4 | 0.80 | 0.42 | 21 |
| | 0.37 | 4.8E-4 | 0.79 | 0.41 | 19 |
| | 0.45 | 5.3E-4 | 0.82 | 0.41 | 19 |
| | 0.48 | 5.4E-4 | 0.82 | 0.41 | 21 |
| 65536 | 0.25 | 8.1E-5 | 0.71 | 0.40 | 37 |
| | 0.28 | 1.1E-4 | 0.74 | 0.41 | 33 |
| | 0.30 | 7.5E-5 | 0.65 | 0.41 | 35 |
| | 0.31 | 6.8E-5 | 0.56 | 0.40 | 35 |
| 1048576 | 0.26 | 3.2E-5 | 0.63 | 0.40 | 65 |
| | 0.32 | 3.0E-5 | 0.56 | 0.40 | 67 |
| | 0.33 | 2.3E-5 | 0.60 | 0.40 | 61 |
| | 0.38 | 4.4E-5 | 0.65 | 0.33 | 103 |

**Table 7:** Characteristics of pdf estimates for the classic Cauchy distribution yielding extreme statistics.

| Sample Size | p-value | KL distance | Figure of Merit | SURD coverage | Lagrange Multipliers |
|---|---|---|---|---|---|
| 256 | 0.48 | 1.7E-2 | 0.87 | 0.48 | 3 |
| | 0.53 | 2.0E-2 | 0.90 | 0.42 | 3 |
| | 0.54 | 1.8E-2 | 0.87 | 0.48 | 3 |
| | 0.54 | 1.2E-2 | 0.88 | 0.50 | 3 |
| 4096 | 0.37 | 7.4E-3 | 0.87 | 0.45 | 4 |
| | 0.38 | 7.4E-3 | 0.87 | 0.45 | 4 |
| | 0.43 | 8.0E-3 | 0.76 | 0.41 | 5 |
| | 0.46 | 7.7E-3 | 0.90 | 0.44 | 4 |
| 65536 | 0.21 | 2.8E-3 | 0.58 | 0.43 | 17 |
| | 0.21 | 2.8E-3 | 0.67 | 0.45 | 17 |
| | 0.26 | 3.1E-3 | 0.71 | 0.41 | 17 |
| | 0.29 | 2.7E-3 | 0.61 | 0.42 | 17 |
| 1048576 | 0.02 | 1.1E-3 | 0.75 | 0.40 | 33 |
| | 0.03 | 1.1E-3 | 0.80 | 0.40 | 35 |
| | 0.03 | 1.1E-3 | 0.76 | 0.40 | 35 |
| | 0.03 | 1.1E-3 | 0.81 | 0.41 | 35 |

**Table 8:** Characteristics of pdf estimates for a discontinuous distribution.

**Figure Captions**

**Figure 1:** QQ plots and SQR plots for SURD of various sample sizes. The left column shows a standard quantile-quantile plot of four independently generated SURD on the interval [0, 1], sorted by position, as a function of their respective population mean positions. From top to bottom, the samples shown are for sizes $N = 2^8, 2^{12}, 2^{16}, 2^{20}$. The right column shows a scaled residual plot of the same SURD generated in the first column across respective rows. The vertical axis is the difference between each sorted value and the expected mean, scaled by the factor $\sqrt{N+2}$ where $N$ is the number of data points in a sample. The SQR plots provide a sample size invariant residual measure that is commensurate with statistical resolution.

**Figure 2:** The universal scoring function is shown as an appropriately shifted probability distribution for SURD log-likelihood based on single order statistics that is empirically calculated. The average log-likelihood for SURD is shown for different sample sizes in different colors as designated by the legend. The scoring function is found to be size invariant shown by the solid black line where distributions for all $N$ are combined into the single curve that is used in this work.

**Figure 3:** A schematic flow chart showing how a random sample of data is iteratively transformed into SURD through hierarchical data partitioning that naturally allows large-scale features to be determined before fine-scale features.

**Figure 4:** Correlation between the Kolmogorov-Smirnov p-value and figure of merit (FOM) for each pdf estimate made per test sample per distribution, yielding (16x8 = 128) points plotted.

**Figure 5:** Assesment of results for the uniform distribution. In the left column, the estimated pdf based on four different samples (shown in green, red, blue and magenta colors) are shown for different sample sizes that range from $N = 2^8, 2^{12}, 2^{16}, 2^{20}$ when moving from the top to bottom rows. These estimates are compared with the true population pdf shown as a black line. The right column shows corresponding SQR-plots across respective rows.

**Figure 6:** Assesment of results for the Lapace distribution. The data format and coloring is the same as Figure 5.

**Figure 7:** Assesment of results for the gamma distribution. The data format and coloring is the same as Figure 5.

**Figure 8:** Assesment of results for the sum of two Gaussian distributions. The data format and coloring is the same as Figure 5.

**Figure 9:** Assesment of results for the 5-fingers perturbing the uniform distribution with a 0.5 weight. The data format and coloring is the same as Figure 5.

**Figure 10:** Assesment of results for the 5-fingers perturbing the uniform distribution with a 0.2 weight. The data format and coloring is the same as Figure 5.

**Figure 11:** Assesment of results for the Cauchy distribution. The data format and coloring is the same as Figure 5.

**Figure 12:** Assesment of results for the discontinuous distribution. The data format and coloring is the same as Figure 5.

**Figure 13:** Log-log plot of the average CPU time to calculate a model pdf as a function of sample size for all eight distributions considered. The averaging was over 100 distinct solutions per sample size.

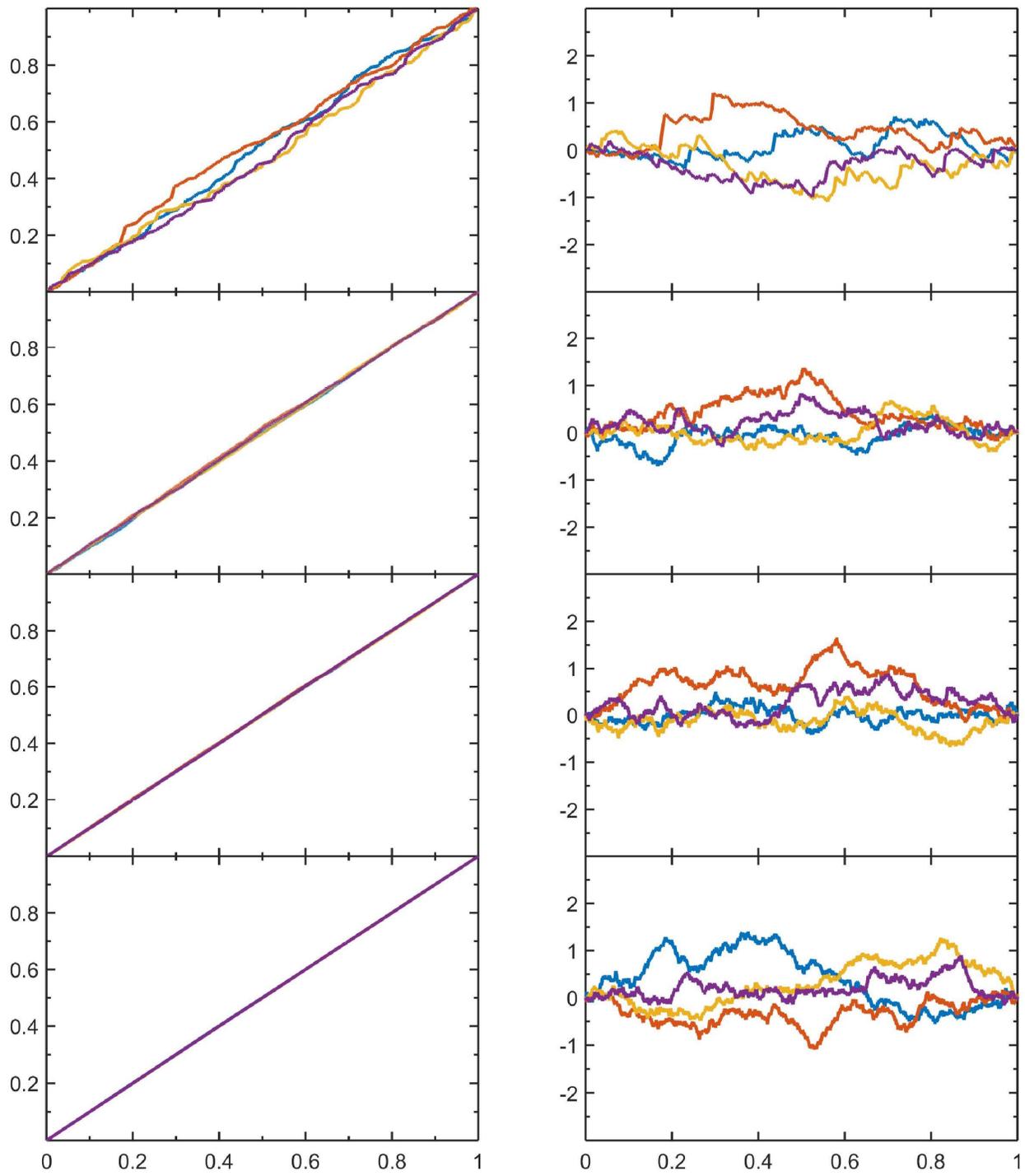

**Figure 1**

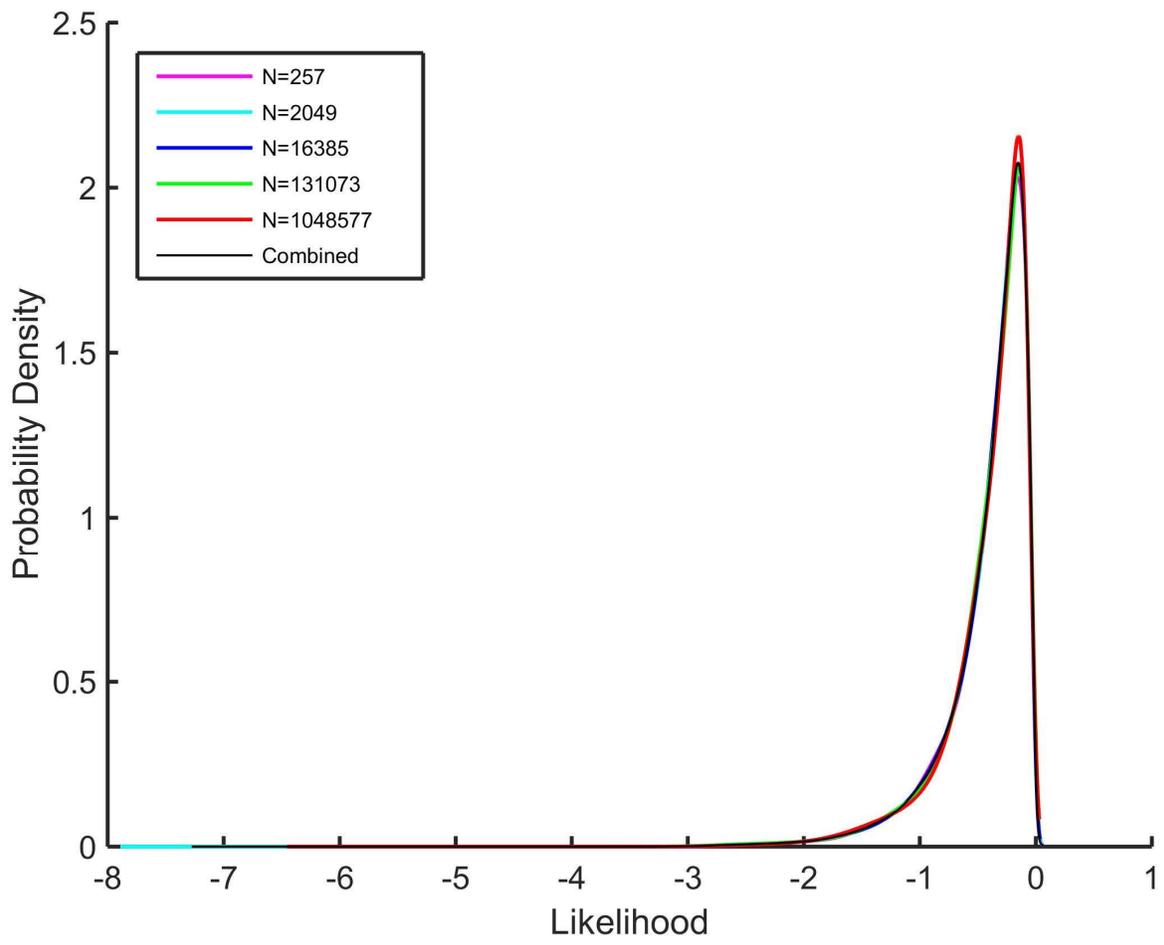

**Figure 2**

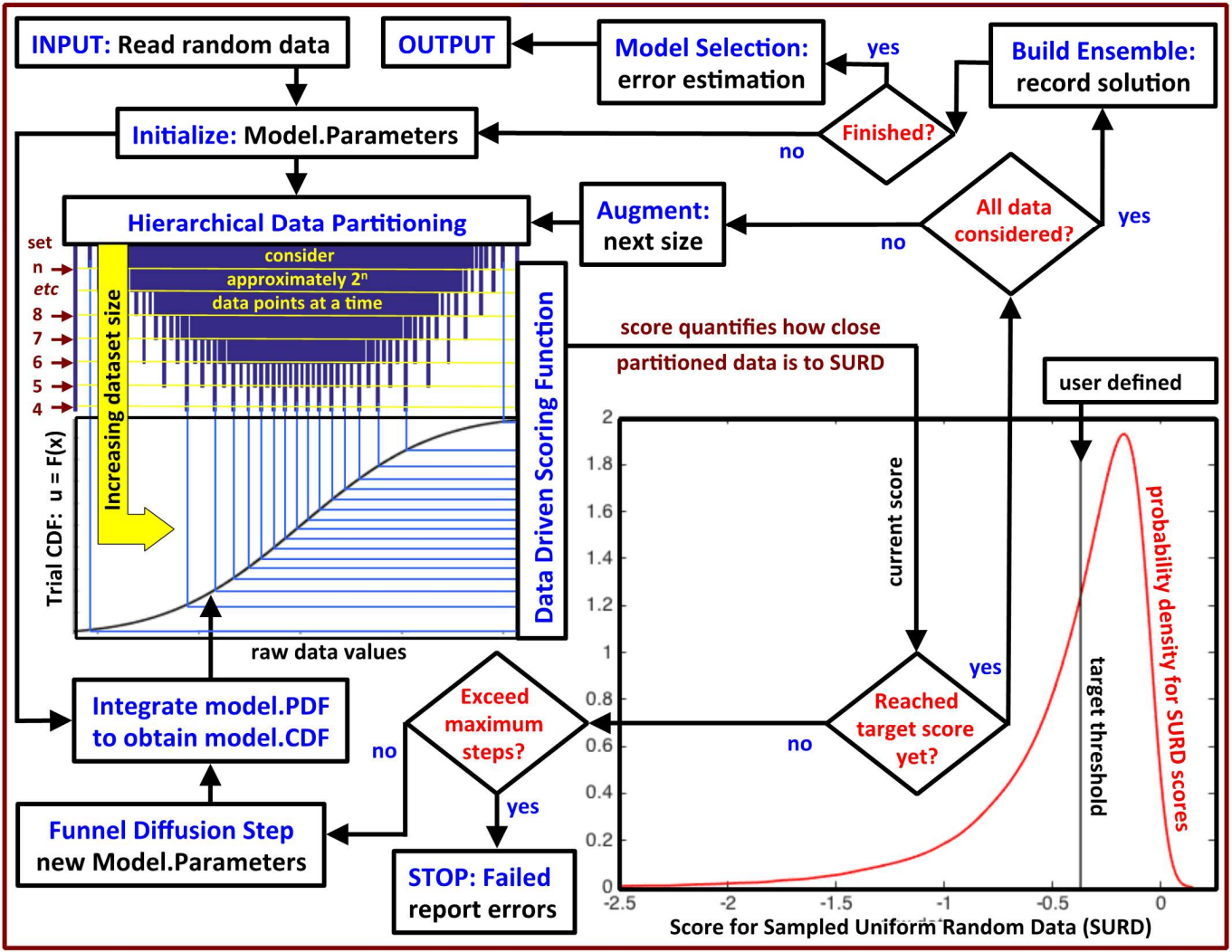

**Figure 3**

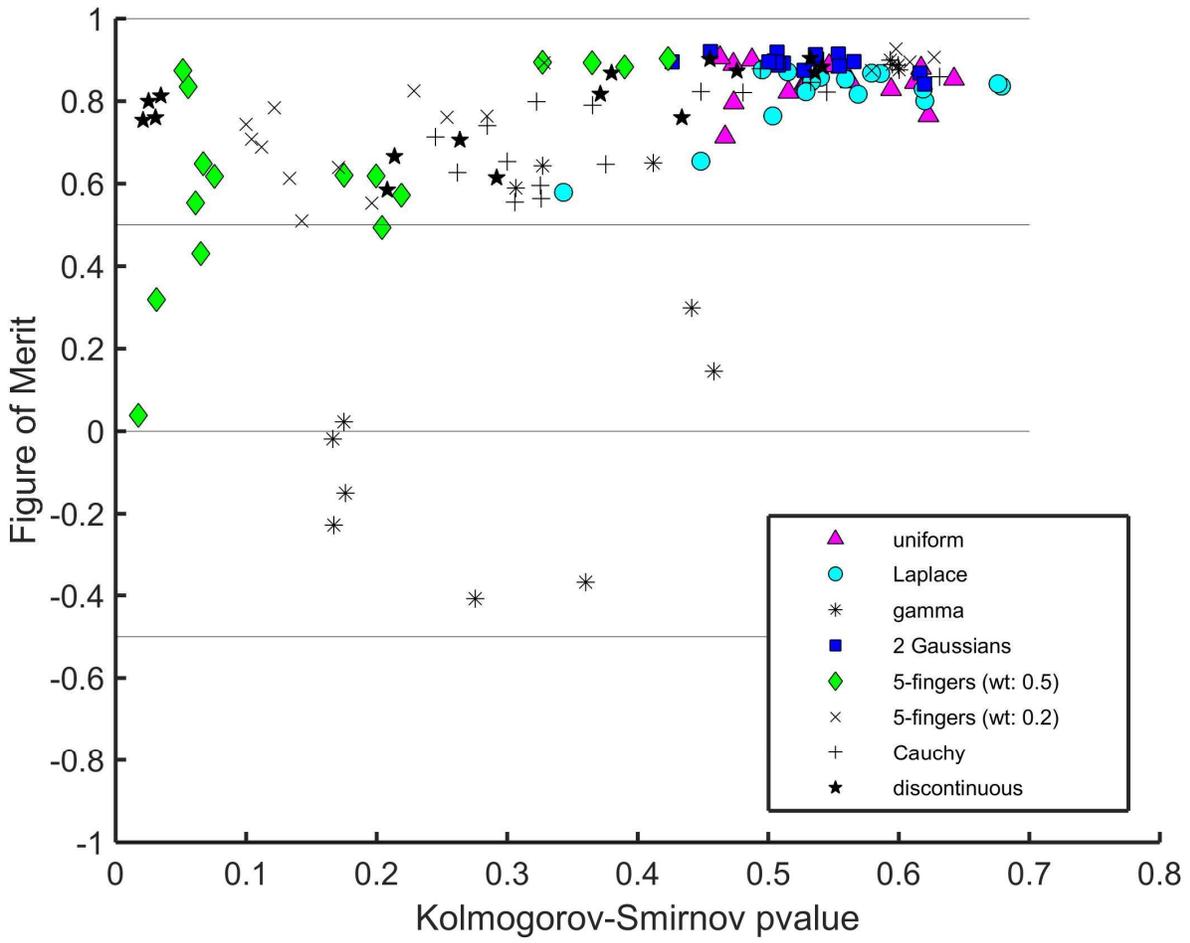

**Figure 4**

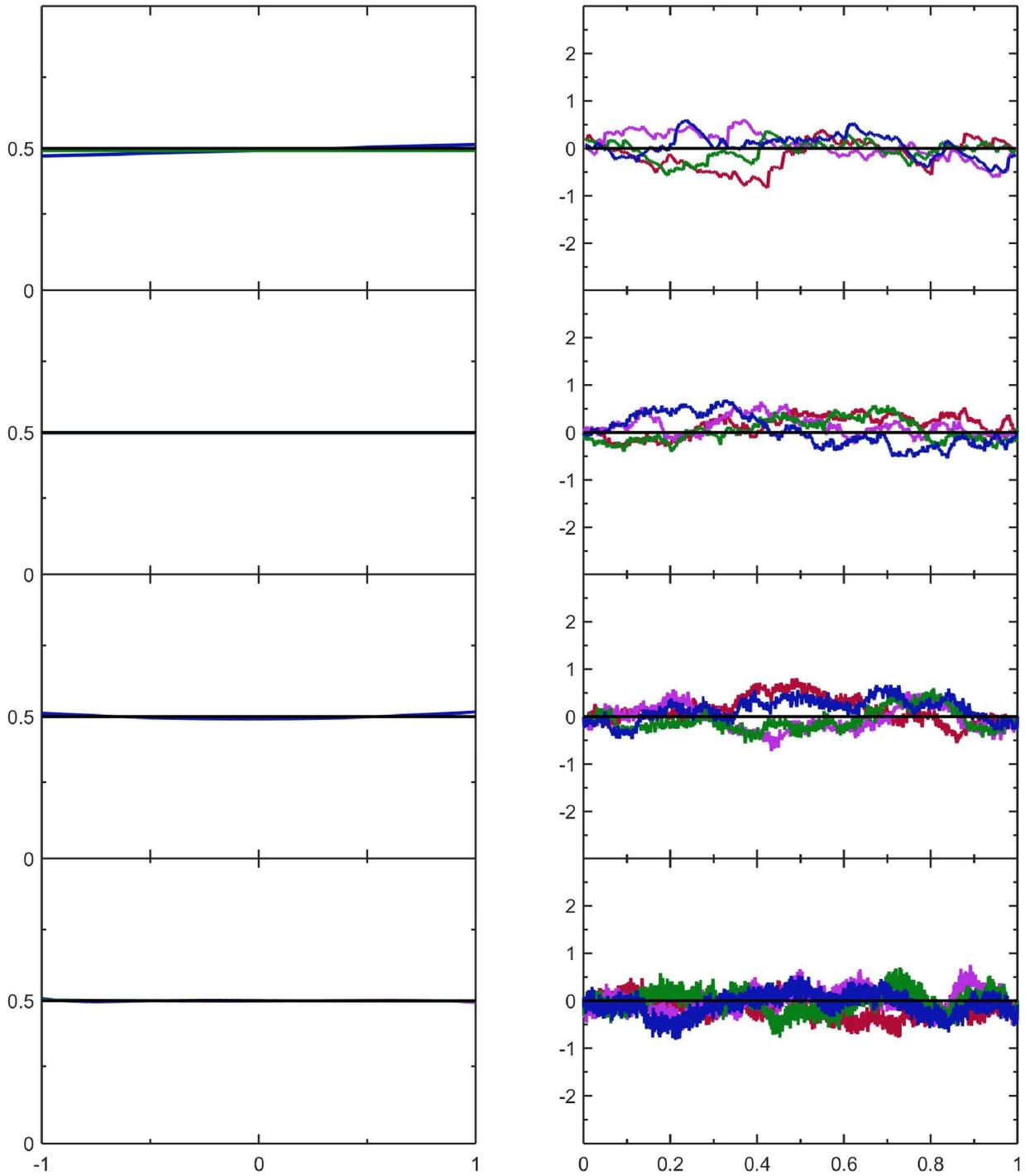

**Figure 5**

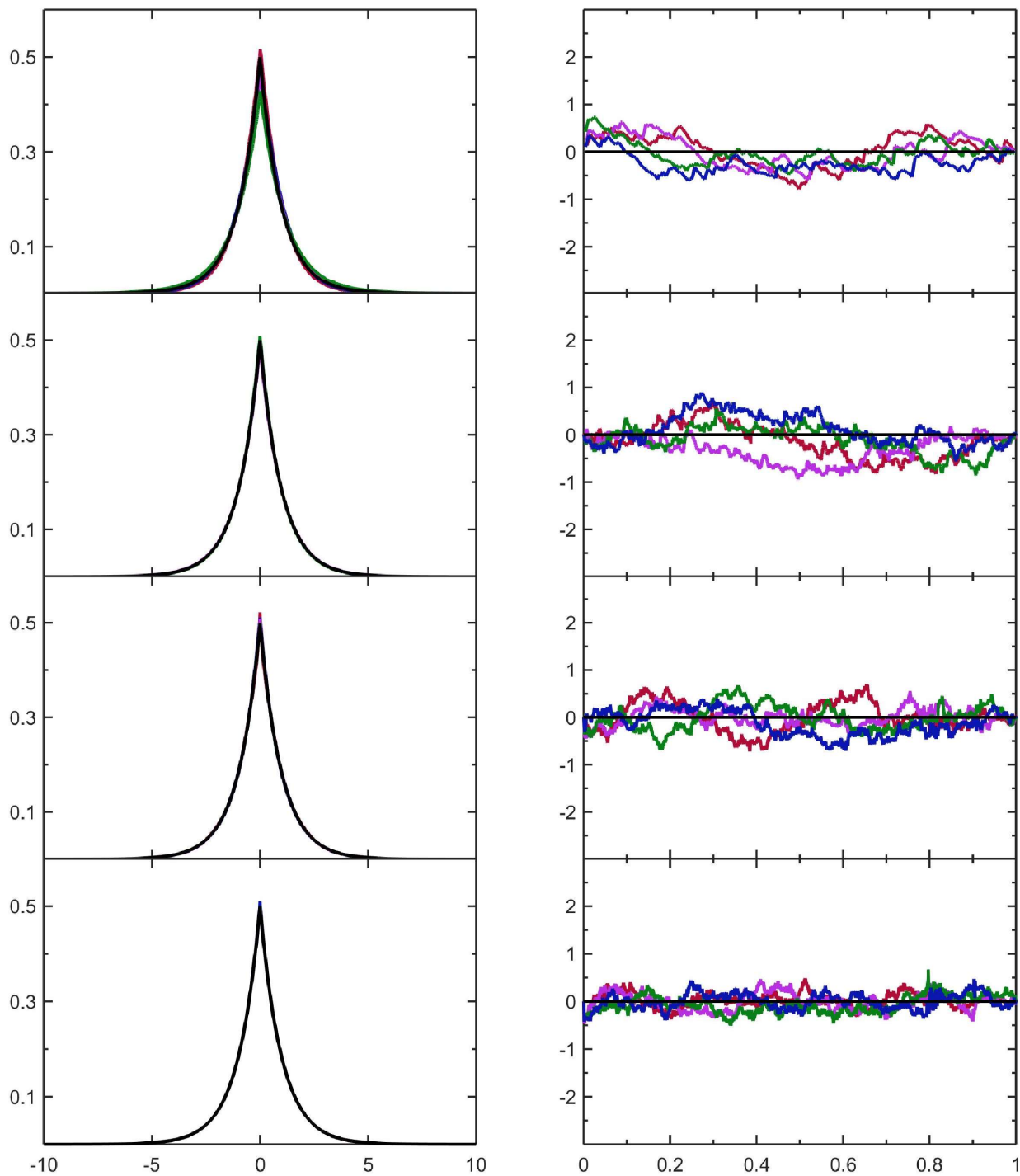

**Figure 6**

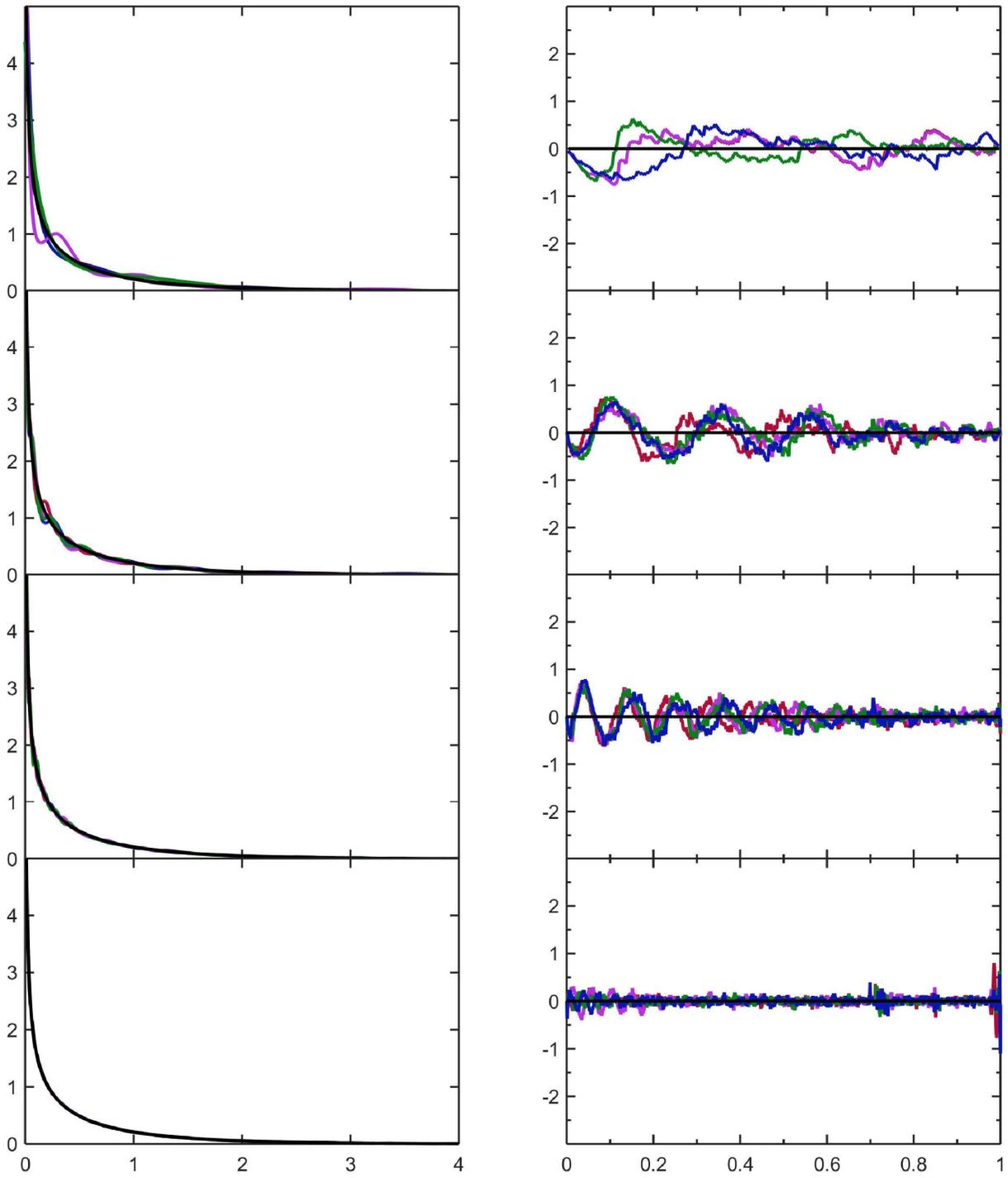

**Figure 7**

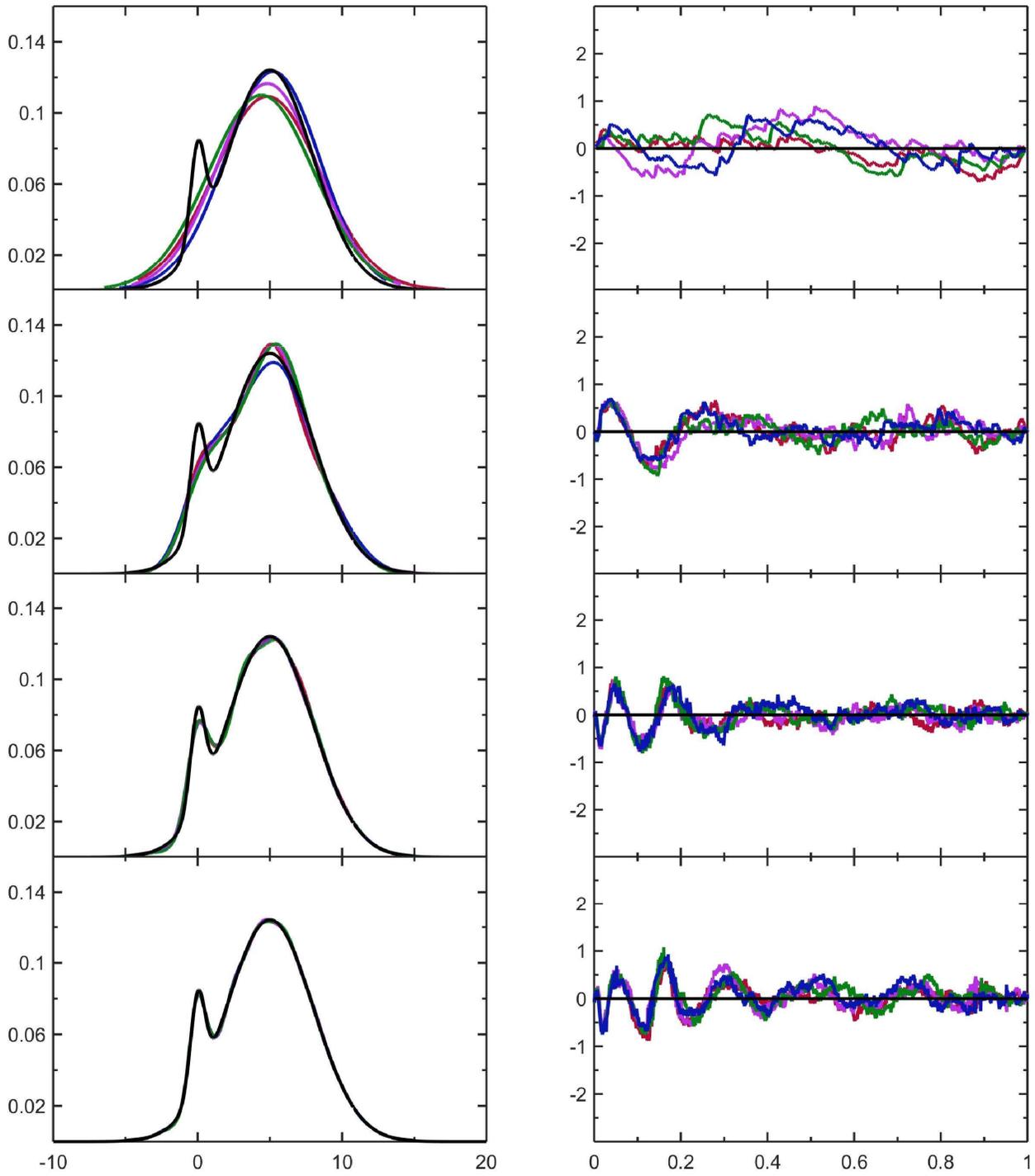

**Figure 8**

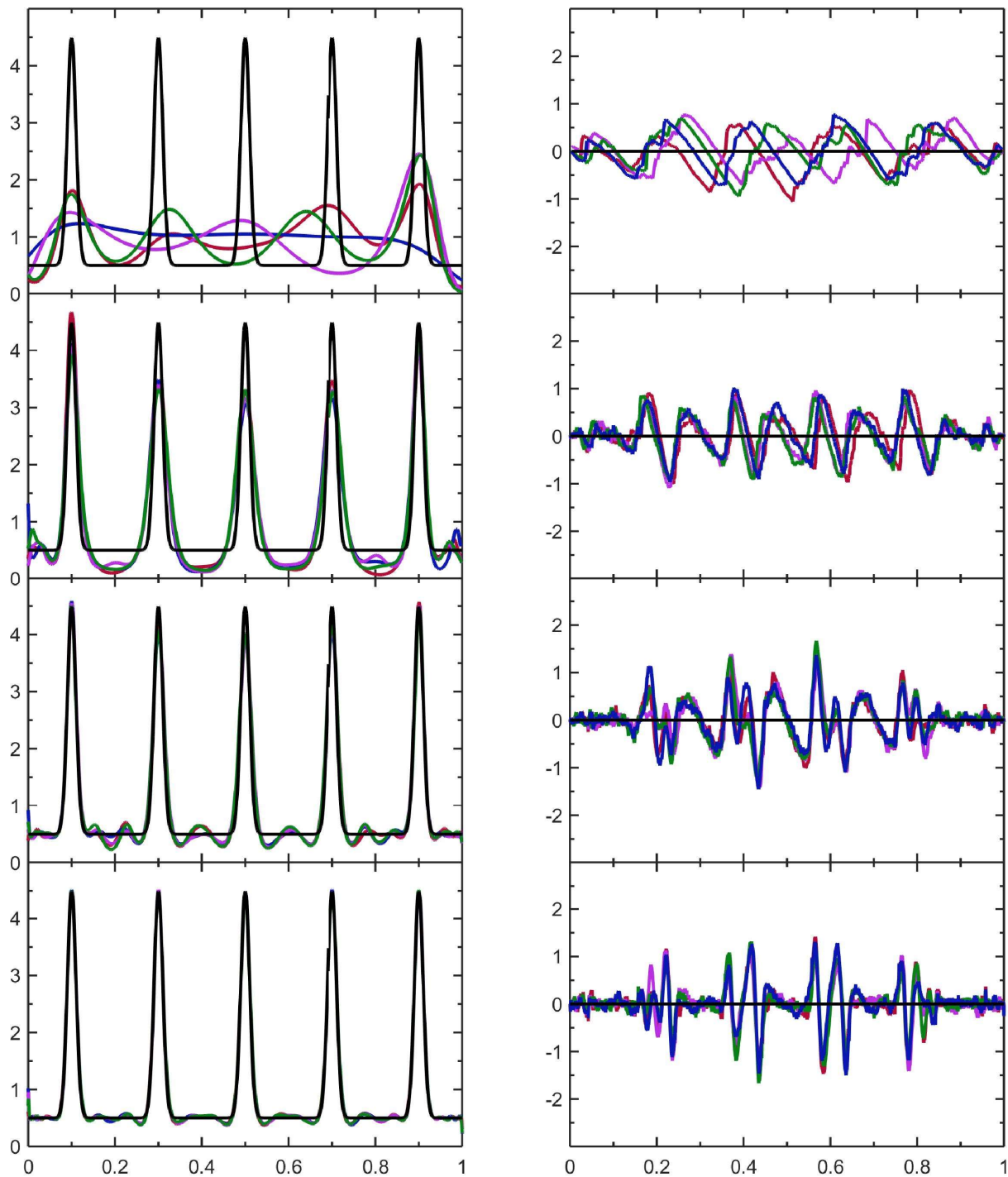

**Figure 9**

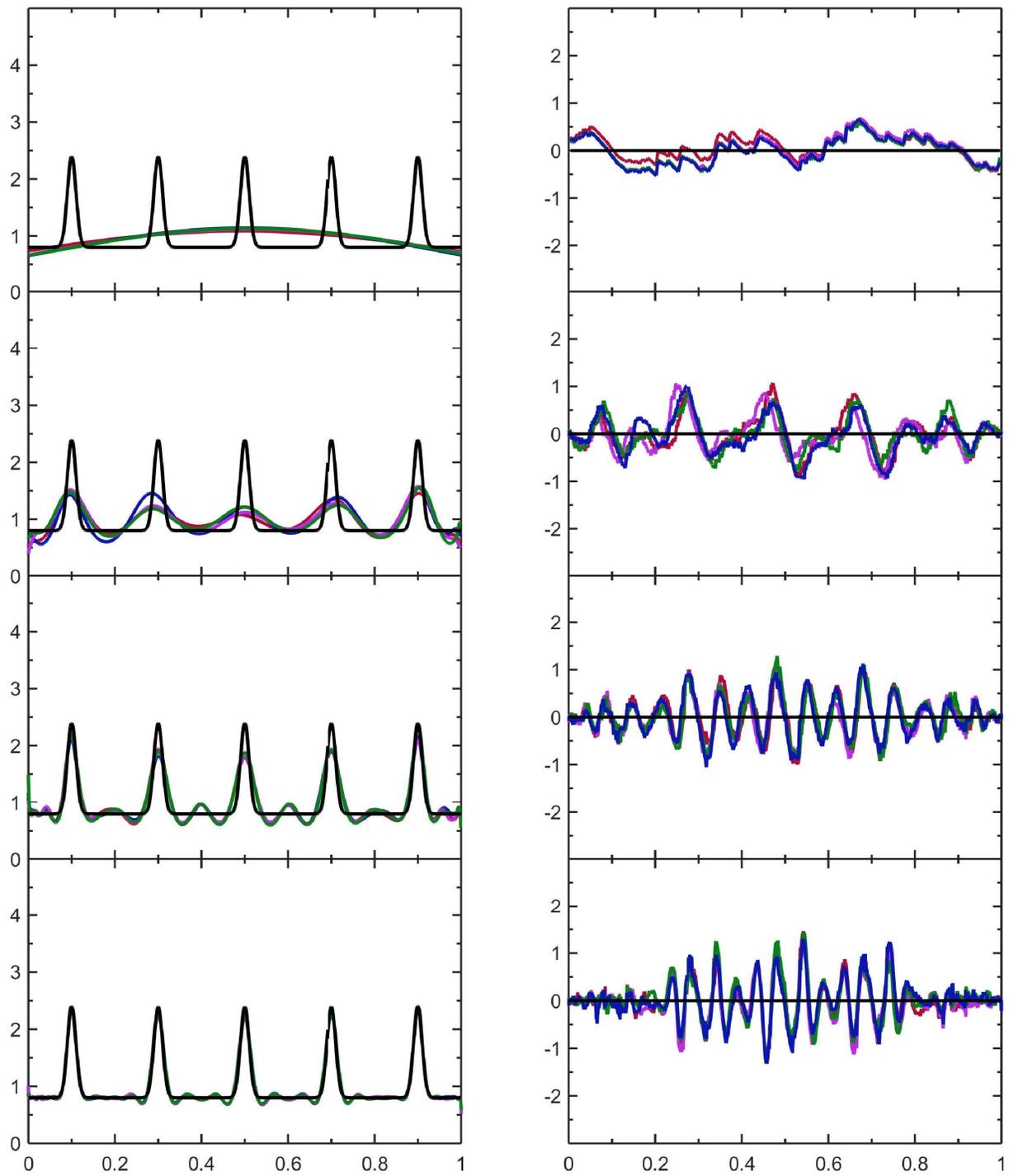

**Figure 10**

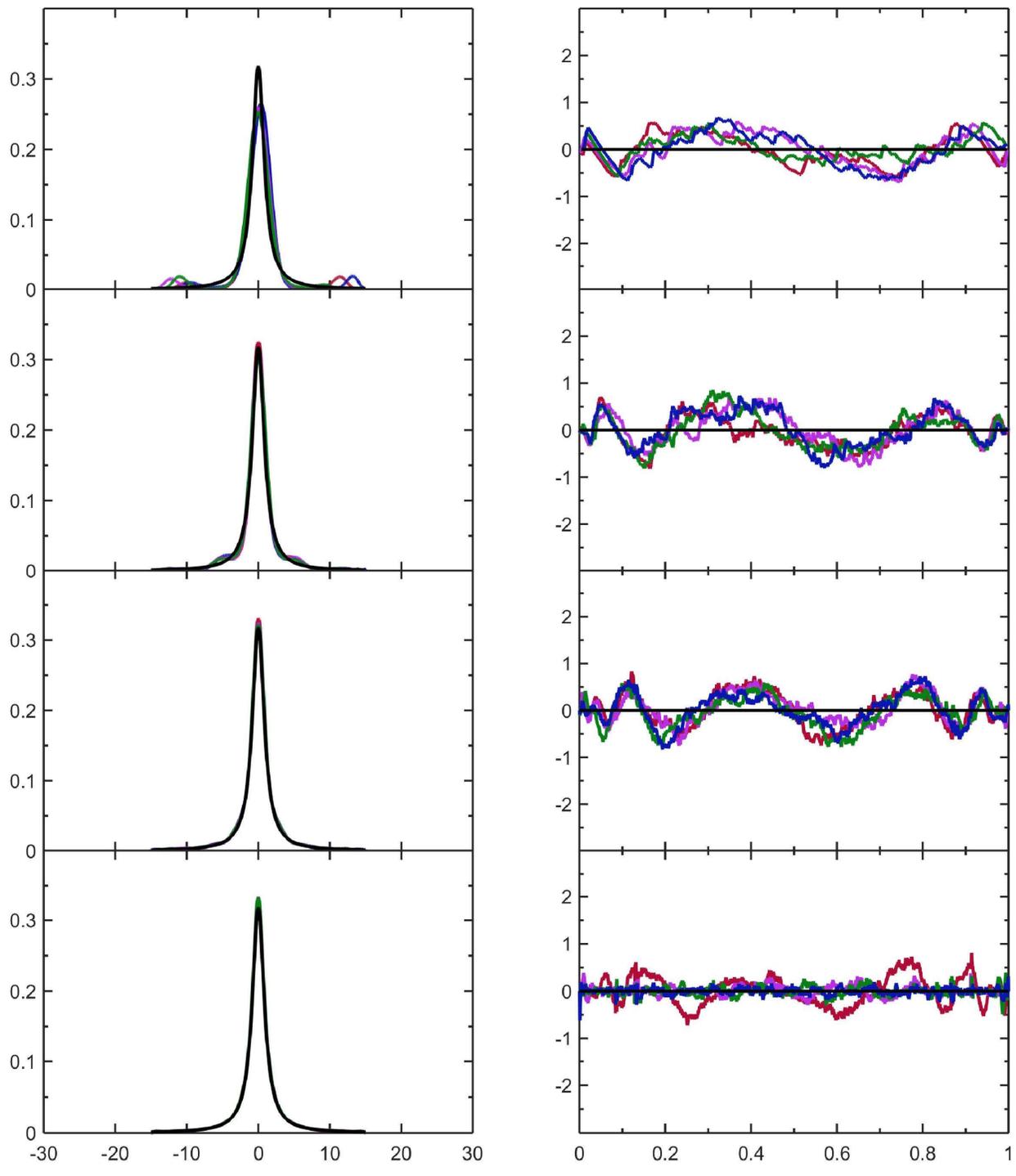

**Figure 11**

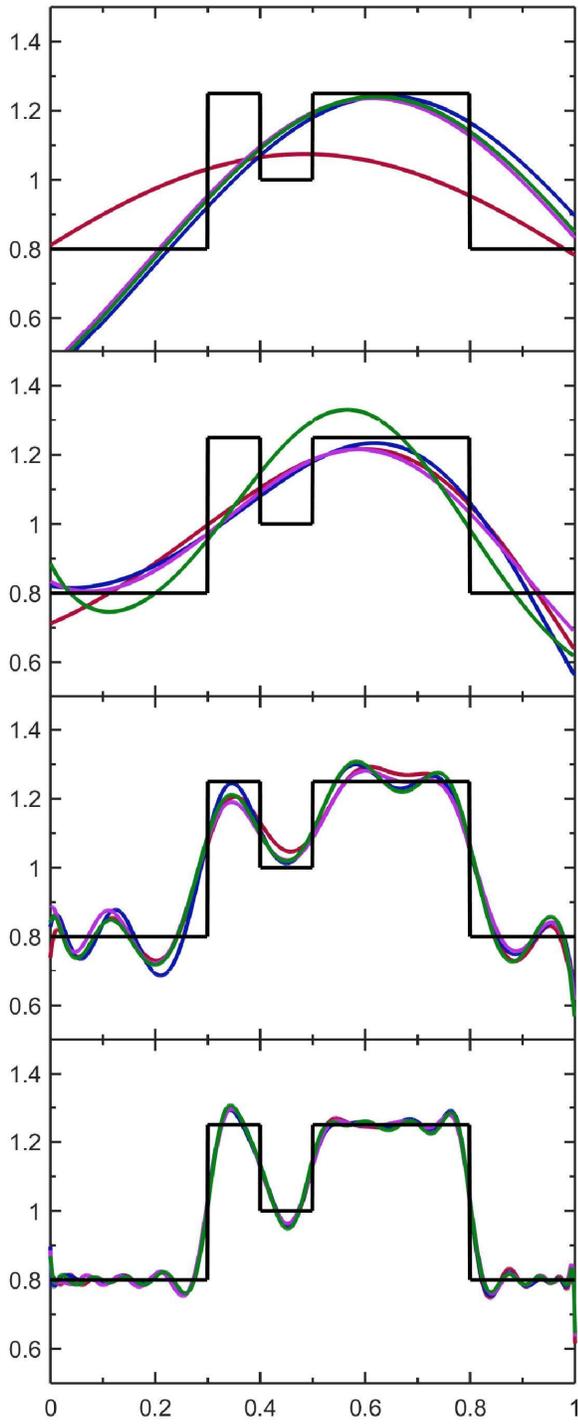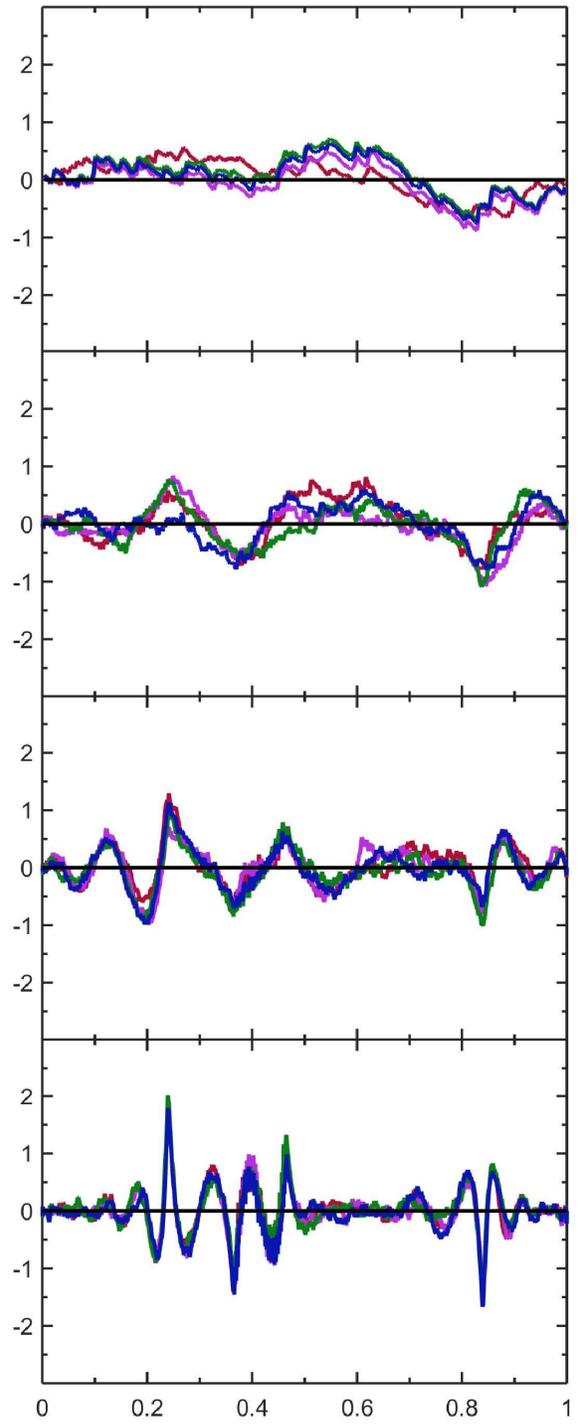

**Figure 12**

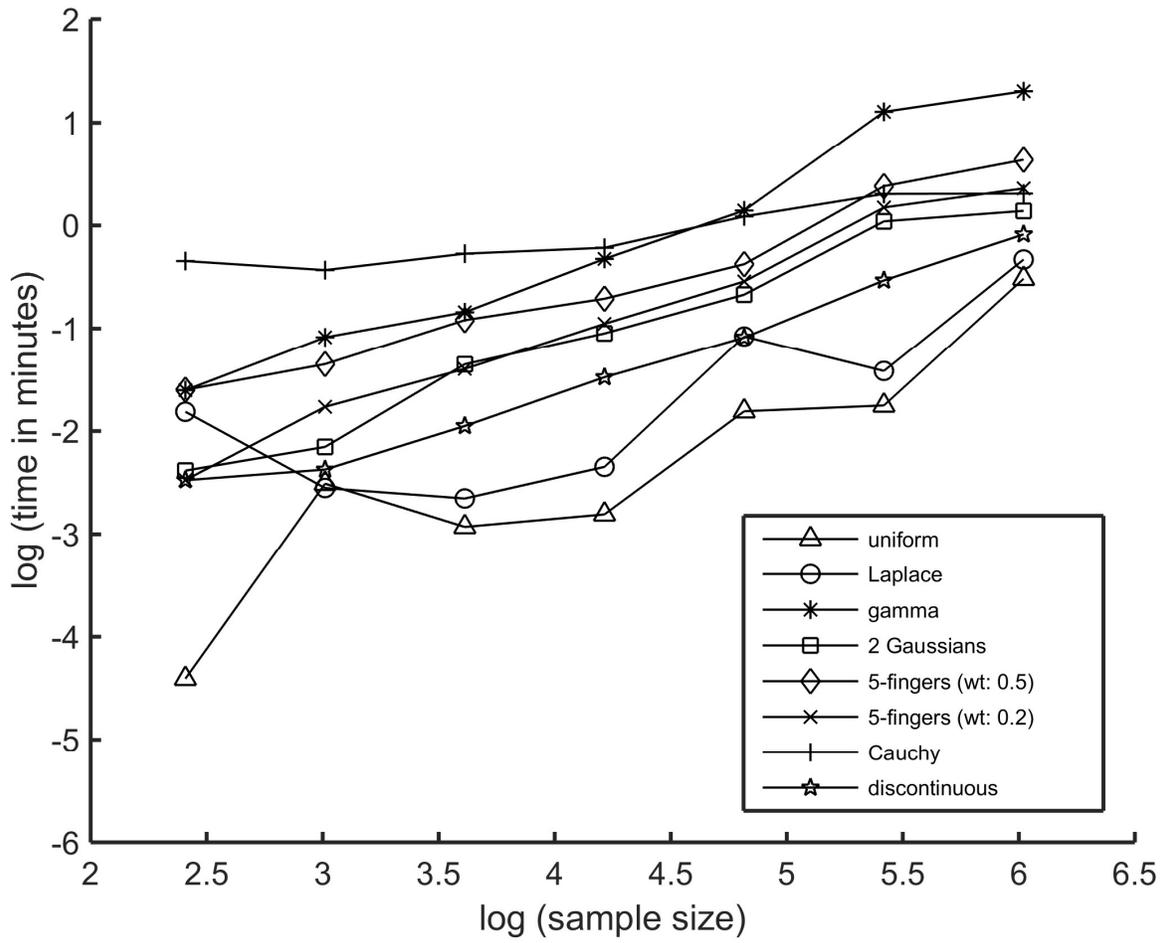

**Figure 13**